\documentclass[onecolumn]{autart}    % Enable this line and disable the 
\usepackage{graphicx}          
\graphicspath{ {./Figures/}}

\usepackage{epstopdf}

\usepackage{algorithm}
\usepackage{algorithmicx}
\usepackage[noend]{algpseudocode}

\usepackage[english]{babel}
\usepackage[utf8]{inputenc}
\usepackage[inline]{enumitem}
\usepackage{ifamath}
\usepackage{tikz}
\usepackage{siunitx}
\usepackage{multicol}
\usepackage{tabu,colortbl}
\usepackage{booktabs}
\usepackage{subcaption}
\usepackage{pgfplots}
\usepackage{wrapfig}
\usepackage{tikz}
\usepackage{tikzscale}
\usetikzlibrary{calc,decorations.pathreplacing,decorations.text,positioning,patterns,topaths,external,3d,spy,matrix,backgrounds}
\usepgfplotslibrary{groupplots}
\pgfplotsset{compat=1.10}
\pgfsetlayers{background,main}
\NewDocumentCommand{\fhighlight}{O{blue!40} m m}{%
  \draw[draw,fill=#1] (#2.north west)rectangle (#3.south east);
}
\newcommand{\includetikzfigure}[2][]{%
  \tikzsetnextfilename{#2}
  \ifstrempty{#1}{
      \includegraphics{#2.tikz}
    }{
      \includegraphics[#1]{#2.tikz}
    }
}

\definecolor{blue}{rgb}{0.25,0.25,0.75}
\colorlet{someblue}{blue!70}
\colorlet{lightblue}{blue!55}
\colorlet{verylightblue}{blue!30}
\definecolor{greyblue}{rgb}{0.36,0.36,0.5}
\colorlet{lightgreyblue}{greyblue!55}
\colorlet{red}{red!90!black}
\colorlet{lightred}{red!65}
\colorlet{verylightred}{red!40}
\colorlet{green}{green!80!black}
\colorlet{lightgreen}{green!65}
\colorlet{verylightgreen}{green!20}
\colorlet{lightgray}{gray!65}

\usepackage{dsfont}
\usepackage[]{units}
\usepackage{xfrac}

\newcommand{\MCP}{\setbf{P}\!\minus\mathbf{MCP}}

\newcommand{\SP}{\setbf{P}\!\minus\mathbf{SP[\mathbb{K}]}}
\newcommand{\PSP}[1]{#1\!\minus\mathbf{SP[\mathbb{K}]}}

\makeatletter
\renewcommand*{\@fnstar}[1]{\ensuremath{\ifcase#1\or \star\or \dagger\or \ddagger\or
    \mathsection\or \mathparagraph\or \|\or **\or \dagger\dagger
    \or \ddagger\ddagger \else\@ctrerr\fi}}
\makeatother

\usepackage[%
  bookmarks=true,
  colorlinks=true,
  pagebackref=false,
  breaklinks=true,
  raiselinks=true,
  pdftitle={Exploiting structure of chance constrained programs via submodularity},
  pdfauthor={Damian Frick, Pier Giuseppe Sessa, Tony A. Wood, Maryam Kamgarpour},
  pdfcreator={pdfLaTeX},
  pdfsubject={},
  pdfkeywords={optimization under uncertainties; randomized methods; chance constrained optimization; submodular minimization; mixed-integer programming}]{hyperref}

\begin{document}

\begin{frontmatter}

\title{Exploiting structure of chance constrained programs via submodularity} % Title, preferably not more 

\thanks[contribinfo]{These authors contributed equally to this work.}
\thanks[grantinfo]{The work of M. Kamgarpour is gratefully supported by Swiss National Science Foundation, under the grant SNSF 200021\_172782.}

\author[Zurich]{Damian Frick\thanksref{contribinfo}}\ead{dafrick@control.ee.ethz.ch},    % Add the 
\author[Zurich]{Pier Giuseppe Sessa\thanksref{contribinfo}}\ead{sessap@control.ee.ethz.ch},               % e-mail address 
\author[Zurich]{Tony A. Wood}\ead{woodt@control.ee.ethz.ch},  % (ead) as shown
\author[Zurich]{Maryam Kamgarpour\thanksref{grantinfo}}\ead{mkamgar@control.ee.ethz.ch}
\address[Zurich]{Automatic Control Laboratory, ETH Zurich, Switzerland}  % Please supply                                              

\begin{keyword}                           % Five to ten keywords,  
    chance constrained optimization; randomized methods; submodular minimization;
\end{keyword}                             % keyword list or with the 
                                          % help of the Automatica 
                                          % keyword wizard

\begin{abstract}                          % Abstract of not more than 200 words.
We introduce a novel approach to reduce the computational effort of solving mixed-integer convex chance constrained programs through the scenario approach. Instead of reducing the number of required scenarios, we directly minimize the computational cost of the scenario program.
We exploit the problem structure by efficiently partitioning the constraint function and considering a multiple chance constrained program that gives the same probabilistic guarantees as the original single chance constrained problem.
We formulate the problem of finding the optimal partition, a partition achieving the lowest computational cost, as an optimization problem with nonlinear objective and combinatorial constraints.
By using submodularity of the support rank of a set of constraints, we propose a polynomial-time algorithm to find suboptimal solutions to this partitioning problem and we give approximation guarantees for special classes of cost metrics.
We illustrate that the resulting computational cost savings can be arbitrarily large and demonstrate our approach on two case studies from production and multi-agent planning.
\end{abstract}

\end{frontmatter}

\section{Introduction}
Mathematical optimization problems are often subject to uncertain parameters affecting the objective or constraint function.
Chance constrained programs (CCPs) \cite{charnes1959,miller1965,pinter1989} provide a framework for solving such optimization problems.
In this framework, the uncertain constraint function is interpreted probabilistically, as a chance constraint, and is allowed to be violated with a given violation probability \cite{prekopa1995}.
Unfortunately, CCPs can be very difficult to solve because they require full knowledge of the probability distribution of the uncertain parameters and are in general non-convex. 

Sample approximations such as the scenario approach \cite{campi2008,calafiore2010}, where a CCP is approximated by a deterministic optimization problem, are not limited to specific distributions and can be applied to a wide range of problems.
In the scenario approach the constraint function is enforced for a number of realizations of the uncertain parameters. However, when many realizations are required these scenario programs can be computationally prohibitive. Numerous efforts have aimed to lessen the computational burden while maintaining the desired probabilistic guarantees of the original CCP.
Most approaches have sought to tighten the sample bounds, to reduce the number of required realizations. In \cite{schildbach2013}, this is achieved by considering the support rank of the constraints.
In \cite{zhang2015}, the structure in the dependence of the constraints on the uncertain parameters is used to further improve these bounds.
In \cite{alamo2015,calafiore2017}, it is proposed to solve a sequence of reduced scenario programs in the framework of sequential probabilistic validation.
Finally, in \cite{margellos2014} the support of the uncertainty is estimated from samples by solving a scenario program of reduced complexity. A robust problem is then solved using linear programming.

A \emph{multiple chance constrained program} (MCP), defined in \cite{schildbach2013}, was constructed in \cite{kariotoglou2016} to minimize the number of constraints of the resulting scenario program.
In this work, we focus on mixed-integer convex CCPs and further develop the perspective of \cite{kariotoglou2016}. We construct multiple chance constrained formulations of the original CCP such that the resulting scenario program has reduced computational cost.
The contributions are as follows:
\begin{enumerate}[label=\arabic*),wide,nosep]
  \item We propose to optimize the partitioning of the constraints into multiple chance constraints, with the objective of minimizing the computational cost of the resulting scenario program.
  Differently from \cite{kariotoglou2016}, where a fixed partition was considered, we also optimize over the partitioning, allowing to more effectively exploit the problem structure. Moreover, we explicitly consider minimizing the computational cost of the resulting scenario program, according to a given cost metric, e.g., the number of floating point operations required for evaluating the constraint. We show that our approach can result in an arbitrarily large reduction of the computational effort.
  \item We show that the \emph{support rank} of a convex constraint is a \emph{monotone submodular} function of the row indices of the constraint function. This implies that the partitioning problem can be reformulated as a minimization of a weakly submodular function (see \refdef{defn:submodularity_ratio} below). We provide a polynomial-time algorithm that exploits this structure and returns suboptimal solutions. We give approximation guarantees for the case where this function is submodular.
\end{enumerate}

\paragraph*{Notation} We denote by $\mathds{1}$ the vector of all ones of appropriate length.
We denote by $\posreals{}$ the non-negative orthant.
We use the Big-O $\bigo(\cdot)$ and Big-Omega $\bigomega(\cdot)$ notation according to Knuth \cite{knuth1976}.
By $\uniform \set{C}$, we denote the continuous (or discrete) uniform distribution supported on the set $\set{C}$. 
Given $x\in \reals{}$, $\ln(x)$ is the natural logarithm of $x$ and $e$ is Euler's number.
Given a finite set $\mathcal{J}$, we denote by $\vert \mathcal{J} \vert$ its cardinality and $\powset{\set{J}}$ the set of all its subsets.
A set function $\gamma : \powset{\set{J}} \rightarrow \reals{}$ is \emph{monotone}, if for every $\set{A},\set{B}$, with $\set{A} \subseteq \set{B} \subseteq \set{J}$ it holds that $\gamma(\set{A}) \leq \gamma(\set{B})$.
A set function $\gamma$ is \emph{submodular} if for every $\set{A},\set{B}$ with $\set{A} \subseteq \set{B} \subset \set{J}$ and every $j \in \set{J} \setminus \set{B}$ it holds that $\gamma(\set{A} \cup \{j\}) - \gamma(\set{A}) \geq \gamma(\set{B} \cup \{j\}) - \gamma(\set{B})$. If this holds with equality, then $\gamma$ is also \emph{modular}.
For simplicity of notation, we define  $\gamma(\set{B} \sep{} \set{A}) := \gamma (\set{A}\cup \set{B}) - \gamma (\set{A})$, and we will use the shorthand $j$, instead of $\{j\}$, if it is clear from context.

\begin{definition}[$\mu$-weak submodularity \protect{\cite{das2011}}]\label{defn:submodularity_ratio}
A monotone set function $\gamma : 2^\mathcal{J} \rightarrow \posreals{}$ is \emph{$\mu$-weakly submodular}, with \emph{submodularity ratio} $\mu \in (0,1]$, if for all $\set{A},\set{B}\subseteq \mathcal{J}$ such that $\set{A} \cap \set{B} = \emptyset$, it holds
\begin{equation*}
  \sum_{j\in \set{B}} \big(\gamma(\set{A} \cup \{j\}) - \gamma(\set{A})\big) \geq \mu \, \big(\gamma(\set{A} \cup \set{B}) - \gamma(\set{A})\big) \,.
\end{equation*}
\end{definition}
\vspace*{-\parskip}
Note: $\gamma$ is $1$-weakly submodular if and only if it is submodular \cite{das2011}. 

\section{Partitioning for the scenario approach}
\label{ccp:sec:section_2}
We consider the following mixed-integer convex chance constrained program (CCP)
\begin{problem}[Mixed-integer convex CCP] \label{ccp:prob:chance_constrained_program}
  \begin{align*}  \min_{(\mathbf{x},\mathbf{y}) \in \setbf{V}} \: & c ^\top \mathbf{x}  \\
    \suchthat \:& \mathbb{P}[f(\mathbf{x},\mathbf{y},\mathbf{w}) \leq 0] \geq 1- \epsilon \,, \nonumber
  \end{align*}
\end{problem}
where $\setbf{V} \subseteq \reals{n} \times \binaries{b}$ is compact and convex for fixed $\mathbf{y}$, $\mathbf{w}$ is an uncertain parameter defined over a probability space $(\mathcalbf{W}, \mathcal{F}, \mathbb{P})$, and the constraint function $f: \reals{n}  \times \binaries{b} \times \mathcalbf{W} \rightarrow \mathbb{R}^r$ is convex in $\mathbf{x}$ for all $\mathbf{y} \in \setbf{Y} := \{ \mathbf{y} \in \binaries{b} \sep{} \exists \mathbf{x} \in \reals{n} \suchthat (\mathbf{x},\mathbf{y}) \in \setbf{V} \}$ and $\mathbf{w} \in \mathcalbf{W}$. The $j$-th component of $f$ is denoted by $f_j$.
The constraint $f$ needs to be satisfied with probability $1-\epsilon$, with $\epsilon \in (0,1)$, where we have used $\mathbb{P}[f(\mathbf{x},\mathbf{y},\mathbf{w}) \leq 0] := \mathbb{P}\{ \mathbf{w} \in \setbf{W} \sep{} f(\mathbf{x},\mathbf{y},\mathbf{w}) \leq 0\}$ as a shorthand.
Note that generic convex objective functions, e.g., expected or worst-case, can be accommodated via a suitable epigraph formulation \cite{calafiore2010}. 

The \emph{violation probability} of \refprob{ccp:prob:chance_constrained_program} for $\mathbf{x}$ and $\mathbf{y}$ is the probability over $\mathcalbf{W}$ for which $\mathbf{x}$,$\mathbf{y}$ violate the constraints in $f$:
\begin{definition}[Violation probability, \protect{\cite[Def.~1]{campi2008}}]\label{ccp:def:constraintViolationProbability}
For a given $(\mathbf{x},\mathbf{y})$, the violation probability of \refprob{ccp:prob:chance_constrained_program} is
\begin{equation*}
  V(\mathbf{x},\mathbf{y}):= \mathbb{P} \{ \mathbf{w}\in \mathcalbf{W} \sep{} \exists j : f_j(\mathbf{x}, \mathbf{y}, \mathbf{w}) > 0 \}\,.
\end{equation*}
\end{definition}
Thus, feasible solutions of \refprob{ccp:prob:chance_constrained_program} satisfy $V(\mathbf{x},\mathbf{y}) \leq \epsilon$.

\subsection{Constraint partitioning}
Without assumptions on $\mathbb{P}$, the feasible set of \refprob{ccp:prob:chance_constrained_program} is non-convex and cannot be expressed in explicit form.
If realizations, or \emph{samples}, of $\mathbf{w}$ are available, the \emph{scenario approach} \cite{campi2008,calafiore2010} can be employed to find points $\mathbf{x}$ which are feasible for \refprob{ccp:prob:chance_constrained_program} with a given \emph{confidence} $1-\beta$, where $\beta \in (0,1)$. Typically, $\beta$ is chosen to be small, e.g., $\beta = 10^{-3}$. For technical reasons, in this paper, we assume $\beta < 1/e \approx 0.36$.
In the resulting \emph{scenario program} the constraint function $f$ is enforced for a finite number samples  of $\mathbf{w}$.
However, often many samples are needed to achieve high confidence and small violation probability, leading to computationally intensive optimization problems. Instead of reducing the number of samples, in this work, we aim to directly reduce the computational effort of evaluating the sampled constraint function.
To do this, we construct a \emph{multiple chance constrained program} (MCP), whose solution is feasible for \refprob{ccp:prob:chance_constrained_program}. We construct this MPC by efficiently partitioning the $r$ rows of the constraint function $f$, allowing us to exploit the structure of $f$.

Given a partition $\setbf{P} := \{\set{P}_i\}_{i=1}^P$ of size $P=|\setbf{P}|$ of the set of row incides $\set{J}:= \{1,\ldots,r\}$ of $f$, we define an MCP
\begin{align*} \label{ccp:eq:multiple_chance_constrained_program}
  \MCP \left\{ \begin{aligned}
    \min_{(\mathbf{x},\mathbf{y}) \in \setbf{V}}\: & c^\top \mathbf{x}  \\
    \suchthat & \mathbb{P}[f_{\set{P}_i}(\mathbf{x},\mathbf{y},\mathbf{w}) \leq 0] \geq (1- \epsilon_i) \quad \forall i=1,\ldots,P\,,
  \end{aligned}  \right.
\end{align*}
where the function $f$ has been partitioned into $\big(f_{\set{P}_1}(\mathbf{x},\mathbf{y},\mathbf{w})\,, \cdots , f_{\set{P}_P}(\mathbf{x},\mathbf{y},\mathbf{w})\big)$, with $f_{\set{P}_i} : \reals{n} \times \binaries{b} \times \setbf{W} \rightarrow \mathbb{R}^{\vert \set{P}_i \vert}$ for $i=1,\ldots,P$.
Moreover, a different violation parameter $\epsilon_i$ has been associated with each constraint $f_{\set{P}_i}$.
As noted in \cite{kariotoglou2016} for the convex case, if $\greekbf{\epsilon} := (\epsilon_1,\ldots,\epsilon_P)$ in $\MCP$ satisfies $\mathds{1}^\transp \greekbf{\epsilon} \leq \epsilon$, then feasible solutions for $\MCP$ are also feasible for \refprob{ccp:prob:chance_constrained_program}. This also holds for the mixed-integer convex case, as outlined in the proof of \refthm{ccp:sec2:thm1} in \refapp{ccp:apdx:mixed-integer}, page~\pageref{ccp:pf:sampling_thm}. Therefore, we approximate $\MCP$ via the scenario approach to obtain feasible solutions for \refprob{ccp:prob:chance_constrained_program}, with a given confidence $1-\beta$.
Due to the presence of binary decisions $\mathbf{y}$, $\MCP$ is more general than the class of MCPs, without binary variables, considered in \cite{schildbach2013}.
Therefore, in \reflem{ccp:lemma:sampling_lemma} we first extend \cite[Section~4.1]{schildbach2013} to the mixed-integer case.
Then, in \refthm{ccp:sec2:thm1} we show that the resulting scenario program can be used to obtain solutions $(\mathbf{x},\mathbf{y})$ that are feasible for \refprob{ccp:prob:chance_constrained_program} with a given confidence $(1-\beta)$.
We will see that the computational cost of the resulting program depends on the partition choice. Hence, we will optimize the partition to minimize this cost.  

\subsection{The scenario approach}\label{ccp:sec:exscapp}
For a partition $\setbf{P}$, the scenario program of $\MCP$ is a convex program where each constraint $f_{\set{P}_i}$ is sampled $K_i$ times. We denote by $\mathbf{w}^{(i,j)}$ the samples related to $f_{\set{P}_i}$ with $j\in\set{K}_i:=\{1,\dots,K_i\}$. We let $K := \sum_{i=1}^{P} K_i$ be the total number of samples. We further assume that $\{\mathbf{w}^{(i,j)} \sep{} i\in\set{I} , j\in\set{K}_i \}$ is a set of i.i.d. random variables, where $\set{I} := \{1,\ldots,P\}$. Given $\mathbb{K} := (K_1,\ldots,K_P)$, the scenario program of $\MCP$ is:
\begin{align*}
  \SP: \left\{ \begin{aligned}
  \min_{ (\mathbf{x},\mathbf{y}) \in \setbf{V}}\: & c^\top \mathbf{x}  \\
  \suchthat & f_{\set{P}_i}(\mathbf{x},\mathbf{y},\mathbf{w}^{(i,j)}) \leq 0 \quad \forall i\in\set{I}, j\in\set{K}_i\,. \nonumber
  \end{aligned} \right.
\end{align*}
Note that $\SP$ depends on the realizations $\mathbf{w}^{(i,j)}$. We assume $\SP$ is feasible for almost all realizations of the multisample and without loss of generality admits a unique solution $(\mathbf{x}^\star[\mathbb{K}],\mathbf{y}^\star[\mathbb{K}])$ which depends on the multisample $\{ \mathbf{w}^{(i,j)} \sep{} i\in\set{I} , j\in\set{K}_i \}$ with product distribution denoted by $\mathbb{P}^K$ \cite{campi2008,calafiore2010,schildbach2013}.

The violation probability of constraint $f_{\set{P}_i}$ is $V_i(\mathbf{x},\mathbf{y}):= \mathbb{P} \{ \mathbf{w}\in \setbf{W} \:\vert\: \exists j \in \set{P}_i : f_j(\mathbf{x}, \mathbf{y}, \mathbf{w}) > 0 \}$. To ensure $V_i(\mathbf{x}^\star[\mathbb{K}], \mathbf{y}^\star[\mathbb{K}])> \epsilon_i$ with high confidence, the number of samples $K_i$ can be chosen depending on the \emph{support rank} of constraint $f_{\set{P}_i}$, which is defined as follows. 
\begin{definition}[Support rank]\label{ccp:def:support_rank_general} 
  We define the support rank as in \cite[Def.~3.6]{schildbach2013}, with the addition of the binary vector $\mathbf{y}$.
  Given a set $\set{A} \subseteq \set{J}$, the \emph{support rank} of the constraint function $f_{\set{A}}$ is defined as 
  \begin{equation*}
    \rho(\set{A}) := n - \dim(L_{\set{A}})\,,
  \end{equation*}
  where $L_{\set{A}}$ is the maximal subspace of $\reals{n}$ not constrained by the function $g(\mathbf{x},\mathbf{y},\mathbf{w}) =\max_{j\in\set{A}} f_j(\mathbf{x},\mathbf{y},\mathbf{w})$\footnote{Formally defined as the maximal linear subspace of $\set{L}_{\max_{j\in\set{A}} f_j}$, given in the appendix on page \pageref{ccp:pf:subspace}.}.
  The support rank function $\rho : \powset{\set{J}} \rightarrow \naturals{}$ therefore maps a set of row indices of $f$ to a non-negative integer.
\end{definition}
Note that the consideration of integer variables does not alter the inherent properties of the support rank and the corresponding proof of \refthm{ccp:sec2:thm1}, because the treatment of $\mathbf{y}$ in the context of the scenario approach relies on the convex problems obtained by fixing $\mathbf{y}$ and considering all possible (feasible) $\mathbf{y}$, see \cite{esfahani2015}.

When the constrains are linear in $\mathbf{x}$ and there are no binary variables, i.e., when $b=0$, specifically, when $f(\mathbf{x},\mathbf{y},\mathbf{w}) = A(\mathbf{w}) \mathbf{x}$ with $A: \setbf{W} \rightarrow \reals{r \times n}$, the support rank is $\rho(\set{A}) = \dim \linspan \{ [A(\mathbf{w})]_j^\transp \mid  j \in \set{A}, \mathbf{w} \in \setbf{W} \}$ \cite{schildbach2013}, where $[A(\mathbf{w})]_j$ is the $j$-th row of $A(\mathbf{w})$. A straight-forward combination of \cite{esfahani2015} with \cite{schildbach2013}, discussed in \reflem{ccp:lemma:sampling_lemma} in \refapp{ccp:apdx:mixed-integer}, yields that if
\begin{equation} \label{ccp:eqn:implicit}
   K_i \in \naturals{} \suchthat 2^b \hspace*{-0.6em}\sum_{l=0}^{\rho(\set{P}_i)-1} \hspace*{-0.2em} \binom{K_i}{l} \epsilon_i^l \left(1-\epsilon_i\right)^{K_i-l} \leq \beta_i \,,
\end{equation} 
then $\mathbb{P}^K [V_i(\mathbf{x}^\star[\mathbb{K}],\mathbf{y}^\star[\mathbb{K}])] > \epsilon_i ] \leq \beta_i$. Moreover, condition~\eqref{ccp:eqn:implicit} is tight for the class of convex fully-supported problems, see \cite{campi2008}.
Using this result, we provide a slight extension of \cite{kariotoglou2016} to the mixed-integer convex case, showing that for a given $\beta$, solutions to \refprob{ccp:prob:chance_constrained_program} can be obtained by solving $\SP$ with appropriate $\mathbb{K} = (K_1, \ldots, K_P)$, as stated in the following theorem. 
\begin{theorem} \label{ccp:sec2:thm1}
Consider \refprob{ccp:prob:chance_constrained_program} and let $\beta \in (0,1)$.  Let $(\mathbf{x}^\star[\mathbb{K}],\mathbf{y}^\star[\mathbb{K}])$ be the unique optimizer of the scenario program $\SP$. Choose $\greekbf{\epsilon}, \greekbf{\beta}  \in (0,1)^P$ such that $\mathds{1}^\transp \greekbf{\epsilon} \leq \epsilon$ and $\mathds{1}^\transp \greekbf{\beta} \leq \beta$.
If each element $K_i$ of $\mathbb{K}$ satisfies \eqref{ccp:eqn:implicit},
then $\mathbb{P}^K [ V(\mathbf{x}^\star[\mathbb{K}],\mathbf{y}^\star[\mathbb{K}]) > \epsilon ] \leq \beta$.
\end{theorem}
\vspace*{-\parskip}
\begin{pf}
See \refapp{ccp:apdx:mixed-integer}, page~\pageref{ccp:pf:sampling_thm}. \hfill\(\qed\)
\end{pf}

Note that specific problem knowledge can be used to improve the bound \eqref{ccp:eqn:implicit}.
For example, $2^b$ can be replaced by the number of feasible binary configurations. This is because the mixed-integer chance constraints of $\MCP$ can be seen as a union of convex chance constraints, where each element corresponds to a feasible configuration of binary variables, see~\cite[Proof of Thm.~4.1]{esfahani2015}.
This number can be substantially lower than $2^b$, e.g. when $\setbf{V} := \{ (\mathbf{x},\mathbf{y}) \in \reals{n} \times \binaries{b} \sep{} \mathds{1}^\transp \mathbf{y} = 1 \}$ there are only $b$ feasible binary configurations, instead of $2^b$.

\subsection{The partitioning problem}\label{sec:partitioning_problem}
The computational effort of solving $\SP$ depends mainly on the constraints of $\SP$. Different metrics can be used to characterize this computational effort, e.g., by considering the \emph{number of constraints} of $\SP$, as proposed in \cite{kariotoglou2016}. In this work, we argue that the cost of evaluating the constraints should be considered explicitly, e.g., by considering the number of \emph{floating point operations} (FLOPs) \cite[p.~12]{golub2013} required to evaluate the constraints, or if the constraints are linear in $\mathbf{x}$, the \emph{number of non-zero elements} (NNZs) of the matrix encoding the constraints.
We consider computational cost metrics $\nu : 2^{\set{J}} \rightarrow \posreals{}$ that map a set of constraint rows to their computational cost. Modularity naturally captures the additive behavior of the above metrics and we therefore restrict ourselves to monotone modular functions $\nu$ such that $\nu(\{j\})>0$ for any $j \in \set{J}$.
Then, given a partition $\setbf{P}$, we define the computational effort associated with $\SP$ as $\sum_{i=1}^P K_i \cdot \nu(\set{P}_i)$, where the sample sizes $K_i$ are selected according to \refthm{ccp:sec2:thm1} as functions of $\epsilon_i$, $\beta_i$ and $ \rho(\set{P}_i)$.
In some cases, instead of reducing the computational complexity of $\SP$, we want to reduce the number of drawn samples, the \emph{sample complexity}. This can also be encoded in a monotone modular metric by setting $\nu(\set{A}) = 1$ for all $\set{A} \subseteq \set{J}$.

In order to reduce the computational cost of solving $\SP$, similarly to \cite{kariotoglou2016}, we make two simplifications:
\begin{enumerate}[label=(\roman*),wide,nosep]
  \item\label{simplification1} We use the \emph{explicit} sample bound \cite{alamo2015}, with the additional term $2^b$ due to the binary variable $\mathbf{y}$
  \begin{equation}\label{ccp:eq:explicit_lb_extended}
   K_i \geq \frac{e}{e -1} \frac{1}{\epsilon_i}\Big( \ln\left(\frac{2^b}{\beta_i}\right) + \rho(\set{P}_i) - 1 \Big) \,.
  \end{equation}
  When \eqref{ccp:eq:explicit_lb_extended} holds, $K_i$ satisfies \eqref{ccp:eqn:implicit} and $\mathbb{P}^K [ V(\mathbf{x}^\star[\mathbb{K}],\mathbf{y}^\star[\mathbb{K}]) > \epsilon ] \leq \beta$ holds. The computational effort associated with $\SP$ is therefore given by
\begin{equation*}
N(\setbf{P}, \greekbf{\epsilon} ,\greekbf{\beta}) := \sum_{i=1}^P K_i \cdot \nu(\set{P}_i) = \frac{e}{e -1} \sum_{i=1}^P \frac{1}{\epsilon_i}\Big( \ln\left(\frac{2^b}{\beta_i}\right) + \rho(\set{P}_i) - 1 \Big) \cdot \, \nu(\set{P}_i)\,.
\end{equation*}
  \item\label{simplification2} We restrict $\greekbf{\beta}$ to depend only on the partition size $P=|\setbf{P}|$, and not on the elements of the partition $\setbf{P}$. We indicate this dependence by the function $\overbar{\greekbf{\beta}}(P)$. Additionally, we require $\mathds{1}^\transp \overbar{\greekbf{\beta}}(P) \leq \beta$ for any $P \leq r$ in order retain the probabilistic guarantees ensured by \refthm{ccp:sec2:thm1}. For the rest of the paper, we fix $\overbar{\greekbf{\beta}}(P) = \frac{\beta}{P} \mathds{1}$.
\end{enumerate}
Obtaining an explicit expression for $K_i$ in \eqref{ccp:eqn:implicit} is challenging and \eqref{ccp:eq:explicit_lb_extended} is the best known lower bound. Moreover, simplification~\ref{simplification2} is reasonable, since the cost $N(\setbf{P}, \greekbf{\epsilon} ,\greekbf{\beta})$ is more sensitive to changes in $\greekbf{\epsilon}$ than in $\greekbf{\beta}$, as also noted in \cite{kariotoglou2016}.
These simplifications lead to the following problem to minimize the computational cost of $\SP$.
\begin{problem}[Partitioning problem]\label{ccp:problem_optimal_R}
\begin{equation*}
\begin{aligned}
\min_{\setbf{P},\,\greekbf{\epsilon}} \quad & N(\setbf{P}, \greekbf{\epsilon} , \tfrac{\beta}{P} \mathds{1}) \\[-0.3em]
\suchthat \,& \greekbf{\epsilon} >0 \:,\, \mathds{1}^\transp \greekbf{\epsilon} \leq \epsilon \,, \\
& \setbf{P} \textnormal{ is a partition of } \{1,\ldots,r\} \,.
\end{aligned}
\end{equation*}
\end{problem}
Note that this generalizes the complexity minimization problem solved in \cite{kariotoglou2016}, where the computational cost metric is the number of rows, i.e., $\nu(\set{P}_i) := |\set{P}_i|$ and only the partition $\big\{\{1\},\ldots,\{r\}\big\}$ is considered.

We motivate our approach by illustrating in \refexa{ccp:example:benefit}, that using a non-trivial partition can lead to significant computational savings. 
In fact, the ratio between the computational costs of the optimal and the trivial partition can be arbitrarily large, as stated in \refprop{ccp:prop:arbitrarily_large}.
\begin{example}\label{ccp:example:benefit}
Consider an instance of \refprob{ccp:prob:chance_constrained_program} with $b=0$ and
\begin{equation*}
 f(\mathbf{x},\mathbf{y},\mathbf{w}) := \begin{tikzpicture}[baseline=-\the\dimexpr\fontdimen22\textfont2\relax]
   \matrix (m)[matrix of math nodes,left delimiter={[},right delimiter={]}] {
     * & \cdots & * & 0 & \cdots & 0 \\
     \vdots & \vdots & \vdots & \vdots & \vdots & \vdots \\
     * & \cdots & * & 0 & \cdots & 0 \\
     * & \cdots & \cdots & \cdots & \cdots & * \\
   };
   \node[fill=red!30,inner sep=0pt] at (m-2-2.center) {$F_1(\mathbf{w})$};
   \node[fill=green!20,inner sep=0pt] at (m-4-3.east) {$F_2(\mathbf{w})$};
   \begin{pgfonlayer}{background}
     \fhighlight[red!30]{m-1-1}{m-3-3}
     \fhighlight[green!20]{m-4-1}{m-4-6}
   \end{pgfonlayer}

 \end{tikzpicture}
 \mathbf{x}\,,
\end{equation*}
where all elements of $F_1(\mathbf{w}) \in \reals{r-1 \times m}$ and $F_2(\mathbf{w}) \in \reals{1 \times n}$ are non-zero, and $f$ has $r$ rows.
We consider two partitions: the trivial partition $\{\set{J} \}$ and $\setbf{P} := \{\set{P}_1, \set{P}_2\}$ with $\set{P}_1 :=\{1,\ldots, r-1 \}$ and $\set{P}_2 := \{r\}$.
In general, $\rho(\set{P}_1) = m$ and $\rho(\set{P}_2) = n$, and therefore $\rho(\set{J}) = n$.
We choose the number of non-zero elements (NNZs) as the cost metric, thus $\nu(\set{P}_1) = rm-m$, $\nu(\set{P}_2) = n$ and $\nu(\set{J}) = rm-m+n$.
Therefore, the NNZs of the constraints of $\PSP{\{\set{J} \}}$ and $\PSP{\setbf{P}}$ are $K(rm-m+n)$ and $K_1 (r-1)m + K_2 n$, respectively, with $K, K_1$ and $K_2$ chosen according to \eqref{ccp:eq:explicit_lb_extended}. \reftab{ccp:tab:benefit} illustrates that depending on the choice of $n$, $m$,and $r$ either partition $\{\set{J} \}$ leads to a problem with lower NNZs, or $\setbf{P}$ does.
\begin{table}[ht]
  \centering
  \begin{tabu*} to \columnwidth {cccc}
    \toprule
    NNZs for & $\PSP{\{\set{J} \}}$ & & $\PSP{\setbf{P}}$ \\
    \cmidrule{2-2} \cmidrule{4-4}
    $m=r=10$, $n=20$ & \cellcolor{green!20} \num{90200} & $<$ & \num{128270}\\
    $m=10$, $r=n=100$ & \num{3652590} & $>$ & \cellcolor{green!20} \num{1715090}\\
    \bottomrule
  \end{tabu*}
  \caption{Number of non-zero elements (NNZs) for \refexa{ccp:example:benefit}, with $\epsilon = 0.05$, $\beta = 10^{-3}$, and $\greekbf{\epsilon} = (\frac{\epsilon}{2}, \frac{\epsilon}{2})$, $\greekbf{\beta} = (\frac{\beta}{2}, \frac{\beta}{2})$ for $\setbf{P}$.}
  \label{ccp:tab:benefit}
\end{table}
\end{example}
\begin{proposition}\label{ccp:prop:arbitrarily_large}
  Given $\epsilon,\beta \in (0,1)$, for any $M > 0$, there exist instances of \refprob{ccp:prob:chance_constrained_program}, partitions $\setbf{P}$ with parameters $\greekbf{\epsilon}$, $\greekbf{\beta}$, and computational cost metrics $\nu$ such that $N(\{\set{J}\},\epsilon,\beta)\,/\,N(\setbf{P},\greekbf{\epsilon},\greekbf{\beta}) \geq M$.
\end{proposition}
\vspace*{-\parskip}
\begin{pf}
  Consider \refexa{ccp:example:benefit} with $m=1$ and $r-1 = n^2$. For $\setbf{P}$, let $\greekbf{\epsilon} := \frac{\epsilon}{2}\mathds{1}$ and $\greekbf{\beta} := \frac{\beta}{2}\mathds{1}$. We can verify that $N(\{\set{J} \},\epsilon,\beta) \geq \tfrac{e}{e-1}\tfrac{1}{\epsilon}n^2(n-1)= \Omega(n^3)$, whereas $N(\setbf{P},\greekbf{\epsilon},\greekbf{\beta}) = \bigo(n^2)$. Thus, for any $M > 0$ there exists an $n$ such that $N(\{\set{J} \},\epsilon,\beta)/N(\setbf{P},\greekbf{\epsilon},\greekbf{\beta}) \geq M$. \hfill\(\qed\)
\end{pf}
\section{Efficient constraint partitioning}
\refprob{ccp:problem_optimal_R} has a nonlinear objective and combinatorial constraints. In fact, the number of all the possible partitions of $\set{J}$ is the Bell number $B_r = \lceil e^{-1} \sum_{k=1}^{2r} \frac{k^r}{k!}\rceil$ and grows exponentially in the number of constraint rows $r$ \cite{rota1964}.
In this section, we propose a greedy algorithm to solve \refprob{ccp:problem_optimal_R}. This algorithm utilizes a particular feature of the support rank function $\rho$, its submodularity, which we prove in \refthm{ccp:thm:submodularity_support_rank}.
\subsection{Submodularity of the support rank function}

The following lemma is used to prove \refthm{ccp:thm:submodularity_support_rank} and shows two properties of the subspace $L_\set{A}$, given in \refdef{ccp:def:support_rank_general}, the definition of the support rank.
\begin{lemma} \label{ccp:lem:subspace}
  Consider a function $f : \reals{n}  \times \binaries{b} \times \setbf{W} \rightarrow \reals{r}$ finite-valued, and convex in its first argument, and two subsets $\set{A},\set{B}$ of rows of $f$, with $\set{A} \subseteq \set{B} \subseteq \set{J}$. Then
  \begin{enumerate*}[label=\roman*)]
    \item\label{ccp:lem:subspace:mon} $L_{\set{B}} \subseteq L_{\set{A}}$, and
    \item\label{ccp:lem:subspace:subm} if $\set{B} \subset \set{J}$ then $\forall j \in \set{J} \setminus \set{B}$: $L_{\set{B} \cup \{j\}} = L_{\set{B}} \cap L_{\set{A} \cup \{j\}}$.
  \end{enumerate*}
\end{lemma}
\vspace*{-\parskip}
\begin{pf}
  The proof is given in \refapp{ccp:apdx:proofs}, page~\pageref{ccp:pf:subspace}.
\end{pf}
\begin{theorem}\label{ccp:thm:submodularity_support_rank}
The support rank function $\rho : \powset{\set{J}} \rightarrow \naturals{}$ associated with the mixed-integer convex constraint function $f$ of \refprob{ccp:prob:chance_constrained_program} is \emph{monotone submodular}.
\end{theorem}
\vspace*{-\parskip}
\begin{pf}
We consider sets $\set{A} \subseteq \set{B} \subseteq \set{J}$.
Using \reflem{ccp:lem:subspace}\ref{ccp:lem:subspace:mon}, we have $L_{\set{B}} \subseteq L_{\set{A}}$, which implies that $\rho(\set{A}) = n - \dim(L_{\set{A}}) \leq n - \dim(L_{\set{B}}) = \rho(\set{B})$, proving monotonicity of $\rho$.
To show submodularity, consider sets $\set{A},\set{B}$ with $\set{A}\subseteq \set{B}\subset \set{J}$ and $j \in \set{J}\setminus\set{B}$. From \reflem{ccp:lem:subspace}\ref{ccp:lem:subspace:subm} it follows that $\dim(L_{\set{B} \cup \{j\}}) = \dim(L_\set{B} \cap L_{\set{A} \cup \{j\}})$. Using \cite[Thm.~4.8]{nering1970}, we have $\dim(L_{\set{B} \cup \{j\}}) = \dim(L_\set{B}) + \dim(L_{\set{A} \cup \{j\}}) - \dim(L_\set{B} + L_{\set{A} \cup \{j\}})$. Thus, $\dim(L_{\set{B} \cup \{j\}}) \geq \dim(L_\set{B}) + \dim(L_{\set{A} \cup \{j\}}) - \dim(L_{\set{A}})$, where we used $\set{A} \subseteq \set{B}$. In fact, $L_\set{B} + L_{\set{A} \cup \{j\}}$ is a subspace of  $L_{\set{A}}$ by \cite[Thm.~4.3]{nering1970} and  \reflem{ccp:lem:subspace}\ref{ccp:lem:subspace:mon}. From the definition it then follows that $\rho(\set{A}\cup \{j\}) - \rho(\set{A}) \geq   \rho(\set{B}\cup \{j\}) - \rho(\set{B})$.\hfill\(\qed\)
\end{pf}
Note that, the number of decision variables $n(\set{A} )$ that enter $f_{\set{A}}$ upper bounds the support rank $\rho(\set{A})$ and is often used as a simple proxy for $\rho(\set{A})$. It is straightforward to show that $n(\set{A})$ is also monotone submodular and the results in this paper also hold for this case.
\subsection{Solving the partitioning problem}\label{ccp:sec:section_3}

To solve \refprob{ccp:problem_optimal_R}, we first provide an explicit expression of the optimal $\greekbf{\epsilon}$, given a partition $\setbf{P}$.
For a given partition size $P$, we define $\sigma : \powset{\set{J}} \rightarrow \posreals{}$ such that
\begin{equation} \label{sigma}
  \sigma_i := \sigma(\set{P}_i) := \big( \ln\left(\tfrac{2^b P}{\beta}\right) + \rho(\set{P}_i) - 1 \big)\cdot \, \nu(\set{P}_i) \,.
\end{equation}
Note that $\sigma(\set{A}) > 0$ for all $\set{A}\neq \emptyset$, and $\sigma(\emptyset)\geq 0$ since $\beta \leq e^{-1}$.
\begin{proposition}\label{ccp:prop:optimal_epsilon}
  Given a partition $\setbf{P}$, \refprob{ccp:problem_optimal_R} has objective value $\frac{e}{e-1} \frac{1}{\epsilon} \big(\sum_{i=1}^P{\sqrt{\sigma_i}}\big)^2$ and a unique minimizer $\greekbf{\epsilon}^\star(\setbf{P}) \in \reals{P}$ with $\epsilon_i^\star(\setbf{P}) = \epsilon \sqrt{\sigma_i}\big(\sum_{j=1}^P{\sqrt{\sigma_j}}\big)^{-1}$.
\end{proposition}
\vspace*{-\parskip}
\begin{pf}
For a given partition $\setbf{P}$, the objective function of \refprob{ccp:problem_optimal_R} can be expressed as a function of $\greekbf{\epsilon}$: 
\begin{equation}\label{ccp:eq:to_substitute}
  N(\setbf{P}, \greekbf{\epsilon} ,\tfrac{\beta}{P} \mathds{1}) = \frac{e}{e -1} \sum_{i=1}^P \frac{\sigma_i}{\epsilon_i}\,.
\end{equation}
It can be verified that the minimizer $\greekbf{\epsilon}^\star$ of \refprob{ccp:problem_optimal_R} must satisfy $\mathds{1}^\transp \greekbf{\epsilon}^\star = \epsilon$, otherwise the cost can be improved by increasing any of the $\epsilon_i$. Moreover, the unique stationary point of \eqref{ccp:eq:to_substitute} subject to $\greekbf{\epsilon}>0$ and $\mathds{1}^\transp \greekbf{\epsilon} = \epsilon$ satisfies
\begin{equation*}
  \epsilon_i^\star = \epsilon \sqrt{\sigma_i}\Big(\sum_{j=1}^P \sqrt{\sigma_j}\Big)^{-1} \quad \text{for } i\in\set{I}\,.
\end{equation*}
The cost tends to infinity as $\epsilon_i \rightarrow 0$ for any $i\in\set{I}$, thus $\greekbf{\epsilon}^\star$ uniquely minimizes \refprob{ccp:problem_optimal_R}.
Moreover, $\greekbf{\epsilon}^\star$ is a function of $\sigma_i$ and therefore of the partition $\setbf{P}$, i.e., $\greekbf{\epsilon}^\star(\setbf{P})$ is a parametric solution to \refprob{ccp:problem_optimal_R}. The objective value is obtained by substituting $\greekbf{\epsilon}^\star(\setbf{P})$ into \eqref{ccp:eq:to_substitute}. \hfill\(\qed\)
\end{pf}

\refprop{ccp:prop:optimal_epsilon} shows how to optimally select $\greekbf{\epsilon}^\star$ as a function of $\setbf{P}$. Therefore, \refprob{ccp:problem_optimal_R} reduces to the optimal selection of $\setbf{P}$, by solving 
\begin{equation}\label{reduction_of_Prob_2}
  \argmin_{\setbf{P}} \sum_{i=1}^{| \setbf{P}|} \gamma(\set{P}_i) \, , 
\end{equation}
with $\gamma(\set{P}_i) := \sqrt{\sigma(\set{P}_i)}$.
Note that \eqref{reduction_of_Prob_2} is a purely combinatorial, in general NP-hard problem \cite{zhao2005}. However, the function $\gamma$ is $\mu$-weakly submodular, as stated in the following lemma.

\begin{lemma}\label{lem:gamma_weak_submodular}
There exists a $\mu \in (0,1]$ such that the set function $\gamma$ is $\mu$-weakly submodular.
\end{lemma}
\vspace*{-\parskip}
\begin{pf}
Recall the definition of $\sigma$ and that the support rank $\rho$ is monotone, non-negative and submodular by \refthm{ccp:thm:submodularity_support_rank}. 
Hence, $\sigma$ is the product of a monotone submodular function with the monotone modular function $\nu$. Therefore, by \reflem{lemma:product_modular_weak_submodular} in \refapp{ccp:apdx:proofs}, there exist a $\mu \in (0,1]$ such that $\sigma$ is $\mu$-weakly submodular. Moreover, by \refprop{prop:composition_with_concave} in \refapp{ccp:apdx:proofs}, $\gamma$ is $\mu$-weakly submodular, as the square root of a monotone $\mu$-weakly submodular function.  \hfill\(\qed\)
\end{pf}

To find partition $\setbf{P}$, a solution to \refprob{ccp:problem_optimal_R}, we propose a greedy algorithm which finds suboptimal solutions in polynomial-time. The structure of the proposed algorithm is motivated by existing results on approximating a class of submodular minimization problems, named \emph{submodular multiway partition} (\textsc{SubMP}) problems. A \textsc{SubMP} is defined as follows. 
\begin{definition}[\textsc{SubMP} \cite{zhao2005}]\label{def:sub_mp}
  Consider a finite set $\set{J}$, a  submodular function $\gamma : \powset{\set{J}} \rightarrow \posreals{}$ and an integer $P \leq \vert \set{J} \vert$. The \textsc{SubMP} problem is to find a partition $\{\set{P}_1,\ldots, \set{P}_P\}$ of $\set{J}$ that minimizes $\sum_{i=1}^P \gamma(\set{P}_i)$.
\end{definition}

In \cite{zhao2005} a \emph{greedy splitting algorithm} (GSA) was developed for solving \textsc{SubMP} in polynomial time with provable approximation guarantees. In this work, we use the GSA to develop a greedy algorithm for finding solutions to \refprob{ccp:problem_optimal_R}. In fact, we can state the following.

\begin{lemma}\label{ccp:proposition_sub_mp}
When $\gamma$ is submodular, \refprob{ccp:problem_optimal_R} subject to $\vert \setbf{P} \vert = P$, for $P \leq r$, can be formulated as a \textsc{SubMP} with a monotone objective function.
\end{lemma}
\vspace*{-\parskip}
\begin{pf}
According to \refprop{ccp:prop:optimal_epsilon}, \refprob{ccp:problem_optimal_R} reduces to solving $\argmin_{\setbf{P}} \sum_{i=1}^{P}\nolimits \gamma(\set{P}_i)$ subject to partition constraints. Hence, when $\gamma$ is submodular, and $\setbf{P}$ is restricted to have size $P$, \refprob{ccp:problem_optimal_R} reduces to a \textsc{SubMP} (see  \refdef{def:sub_mp}).  Moreover, $\gamma$ is the square root of the product of two monotone non-negative functions and is therefore also monotone.  \hfill\(\qed\)
\end{pf} 

We propose \refalg{ccp:algorithm_paritioning} for finding suboptimal solutions to \refprob{ccp:problem_optimal_R}.
We iterate over the partition size $P$, and for each $P$ we run the GSA algorithm of \cite{zhao2005}. Although developed for submodular objectives, the GSA and its subroutines \cite{queyranne1998} run independently of the submodularity of $\gamma$. Hence, as stated in \refthm{ccp:thm_optimal_partitioning}, \refalg{ccp:algorithm_paritioning} has a run time that is polynomial in $r$, the number of constraint rows. Moreover, when $\gamma$ is submodular, approximation guarantees are available.

\begin{algorithm}[h]
  \caption{Efficient partitioning}\label{ccp:algorithm_paritioning} 
    \begin{algorithmic}[1]
      \Require{Desired parameters $\epsilon, \beta$ and cost function $\gamma$}
    \State Initialize: $N^\star \gets N(\{\set{J}\}, \epsilon,\beta)$, $\setbf{P}^\star \gets \{\set{J}\}$      
       \For{\texttt{$P=2,\ldots, r$}}\label{ccp:alg:line:iterate}
       \State $\setbf{P}$ $\gets$ \textsc{GSA}($P,\gamma$) 
       \State $\overbar{N} \gets N(\setbf{P},\greekbf{\epsilon}^\star(\setbf{P}),\tfrac{\beta}{P} \mathds{1})$ \Comment{compute cost}
         \If{$\overbar{N} < N^\star$}          
           $\: \setbf{P}^\star \gets \setbf{P}$ , $N^\star \gets \overbar{N}$
         \EndIf 
       \EndFor
      \State $\greekbf{\epsilon}^\star \gets \greekbf{\epsilon}^\star(\setbf{P}^\star)$ \Comment{apply \refprop{ccp:prop:optimal_epsilon}}
      \State \Return $\setbf{P}^\star, \greekbf{\epsilon}^\star$
    \end{algorithmic}
\end{algorithm}

\begin{theorem}\label{ccp:thm_optimal_partitioning}
Given $r \geq 2$, \refalg{ccp:algorithm_paritioning} returns $\setbf{P}^\star,\greekbf{\epsilon}^\star$, a feasible solution to \refprob{ccp:problem_optimal_R}, with a number of evaluations of $\gamma$ which is polynomial in $r$. Moreover, if $\gamma$ is submodular, then the solution has a $(2-2/r)^2$-approximation guarantee.
\end{theorem}
\vspace*{-\parskip}
\begin{pf}
\refalg{ccp:algorithm_paritioning} has $r-1$ iterations. At each iteration the GSA finds a partition of size $P$. Regardless of submodularity of $\gamma$, the GSA requires $\bigo(P r^3)$ evaluations of $\gamma$ \cite{zhao2005,queyranne1998}. Moreover, if $\gamma$ is monotone submodular, GSA solves a \textsc{SubMP} returning a $(2-2/P)$-approximation \cite{zhao2005}. Using \reflem{ccp:proposition_sub_mp} and \refprop{ccp:prop:optimal_epsilon}, we therefore obtain that the resulting $\setbf{P}^\star,\greekbf{\epsilon}^\star$ correspond to a $(2-2/r)^2$-approximation in the worst-case. The exponent 2 in the approximation guarantee is due to the fact that \textsc{SubMP} has an objective which is the square root of the original objective, see  \refprop{ccp:prop:optimal_epsilon}. \hfill\(\qed\)
\end{pf}

Note that the guarantees of Theorem \ref{ccp:sec2:thm1} directly apply to $\setbf{P}^\star$ and $\greekbf{\epsilon}^\star$ returned by \refalg{ccp:algorithm_paritioning}. Moreover, if the sample complexity is chosen as the cost metric, $\sigma$ is submodular by submodularity of the support rank $\rho$. Thus $\gamma$ is submodular, being the square root of a submodular function. Then \refalg{ccp:algorithm_paritioning} gives an \emph{$\alpha$-approximation} guarantee with $\alpha = (2-2/r)^2$, i.e., the solution of \refalg{ccp:algorithm_paritioning} has an objective value at most a factor $\alpha$ larger than the optimum.

\section{Case studies} \label{ccp:sec:section_4}
In this section we demonstrate our theoretical results on two numerical case studies by applying the results developed in the previous sections to different instances of these problems.
In both case studies, the number of samples to build the scenario programs is determined using the tighter implicit bound, that is the smallest $K_i$ that satisfy \eqref{ccp:eqn:implicit}, which are found using a binary search.
In the first case study we consider a production planning example, similar to \cite{escudero1993}, without mixed-integer constraints, i.e., $b=0$. 
For different instances of this problem, we compare three computational cost metrics $\nu$. In particular, we show that the number of floating point operations (FLOPs) of $\SP$ can be reduced via partitioning, and this leads to a reduction of the actual computational effort, measured as the median solver time.
In the second case study, we consider a multi-agent planning task involving mixed-integer constraints and illustrate that the theoretical guarantees of \reflem{ccp:lemma:sampling_lemma} and \refthm{ccp:sec2:thm1} hold by examining the empirical violations.

\subsection{Production planning with capacity constraints} \label{ccp:case_study:production_planning}
The goal is to optimally plan the production of $m$ products. To meet the known demands $\mathbf{d} \in \reals{m}$, products can be either produced internally or procured externally. The $n$ machines available for internal production have an uncertain capacity matrix $\mathbf{A} \in \reals{n \times m}$, where $\mathbf{A}_{ij}$, the $i,j$-th element of $\mathbf{A}$, is the capacity of machine $i$ needed to produce product $j$.
The vectors $\mathbf{x}_1, \mathbf{x}_2 \in \reals{m}$ are the amounts of internal production and external procurement of each product, respectively.
The total capacity of each machine is normalized to one, i.e, $\mathbf{A} \mathbf{x}_1 \leq \mathds{1}$.
In addition to $\mathbf{A}$, the production costs $\mathbf{q} \in \reals{m}$ and the procurement prices $\mathbf{p} \in \reals{m}$ per unit product, as well as the unmet demand costs $\mathbf{u} \in \reals{m}$ are uncertain.
The objective is to minimize the worst-case overall cost which can be formulated as the following CCP:
\begin{align}\label{ccp:example:robust_program}
  \min_{\substack{ \tau \in \mathbb{R},\, (\mathbf{x}_1, \mathbf{x}_2)\in \setbf{X} }} \:& \tau  \\
  \suchthat \:& \mathbb{P}\left[\begin{aligned} \mathbf{A} \,\mathbf{x}_1 &\leq \mathds{1}\,,\\
                                            \tau &\geq  l(\mathbf{x}_1,\mathbf{x}_2,\mathbf{q},\mathbf{p},\mathbf{u})
                          \end{aligned}\right] \geq 0.9\,, \nonumber
\end{align}
with objective function $l(\mathbf{x}_1,\mathbf{x}_2,\mathbf{q},\mathbf{p},\mathbf{u}) := \mathbf{q}^\top\mathbf{x}_1 + \mathbf{p}^\top\mathbf{x}_2 +\mathbf{u}^\top(\mathbf{d}-\mathbf{x}_1-\mathbf{x}_2)$ and the deterministic polyhedral set $ \setbf{X}  := \big\{ (\mathbf{x}_1, \mathbf{x}_2) \in \mathbb{R}^{2m} \sep{\big} \mathbf{x}_1,\: \mathbf{x}_2 \geq 0 ,\: \mathbf{x}_1 + \mathbf{x}_2 \leq \mathbf{d} \big\}$.
The uncertain parameters $\mathbf{A},\mathbf{q},\mathbf{p}$, and $\mathbf{u}$ are distributed as follows: $\mathbf{A}_{ij} \sim \uniform [0.0015, 0.0025]$, $q_j \sim \uniform [\overbar{q}_j, \frac{5}{3}\overbar{q}_j]$, $p_j \sim \uniform [\overbar{p}_j, \frac{5}{3}\overbar{p}_j]$ and $u_j \sim \uniform [\overbar{u}_j, \frac{5}{3}\overbar{u}_j]$, where $\overbar{q}_j$, $\overbar{p}_j$ and $\overbar{u}_j$ represent deterministic nominal values.
We assume that the distributions of these uncertain parameters are not known and that only samples from the distributions are available.
A particular random instance of the demands and nominal values has been chosen as $d_j \sim \uniform \{1,\ldots,100\}$, $\overbar{q}_j \sim \uniform [0, 0.3]$, $\overbar{p}_j \sim \uniform [0,0.6]$ and $\overbar{u}_j \sim \uniform [0,0.9]$.

Problem~\eqref{ccp:example:robust_program} has the form of \refprob{ccp:prob:chance_constrained_program}, with $\epsilon=0.1$. We want to obtain a $(1-\beta)$-confidence solution via the scenario approach, with $\beta = 10^{-3}$. We apply the results developed in the previous sections with the aim of reducing the computational burden of the corresponding scenario program.
We consider problem instances for $m \in \{10,20,30,40\}$ products and $n \in \{10,20, \ldots, 300\}$ machines. For each $m$ and $n$, we run \refalg{ccp:algorithm_paritioning} to obtain a partition $\setbf{P}^\star$ and corresponding $\greekbf{\epsilon}^\star$.
For computational reasons we only consider partitions with at most four elements in \refalg{ccp:algorithm_paritioning}, i.e., we iterate for $P=2,3,4$ in line 2.
Note that the number of partitions for $P = 4$ and $n=300$ is in the order of $10^{179}$ which means that an exhaustive search is not tractable even for small $P$.
For different computational cost metrics $\nu$, we compare the cost $N\big(\{\set{J}\},\epsilon,\beta\big)$ of the trivial partition $\{\set{J}\}$ to the cost $N\big(\setbf{P}^\star,\greekbf{\epsilon}^\star, \tfrac{\beta}{|\setbf{P}^\star|}\mathds{1}\big)$. We additionally compare this to the actual computational effort of solving $\PSP{\{\set{J}\}}$ and $\PSP{\setbf{P}^\star}$.
To examine the actual computational effort, we consider $20$ random instances of \eqref{ccp:example:robust_program} for each choice of $m$ and $n$.
For each instance, we solve the scenario programs $\PSP{\{\set{J}\}}$ and $\PSP{\setbf{P}^\star}$ on an Intel~i7 CPU with $4$ cores at \unit[2.8]{GHz} and \unit[8]{GB} of memory, using YALMIP \cite{lofberg2004} and Gurobi \cite{gurobi}.

\begin{figure}[ht]
  \centering
  \begin{subfigure}[b]{0.485\linewidth}
  \centering
    \includegraphics[width=\linewidth]{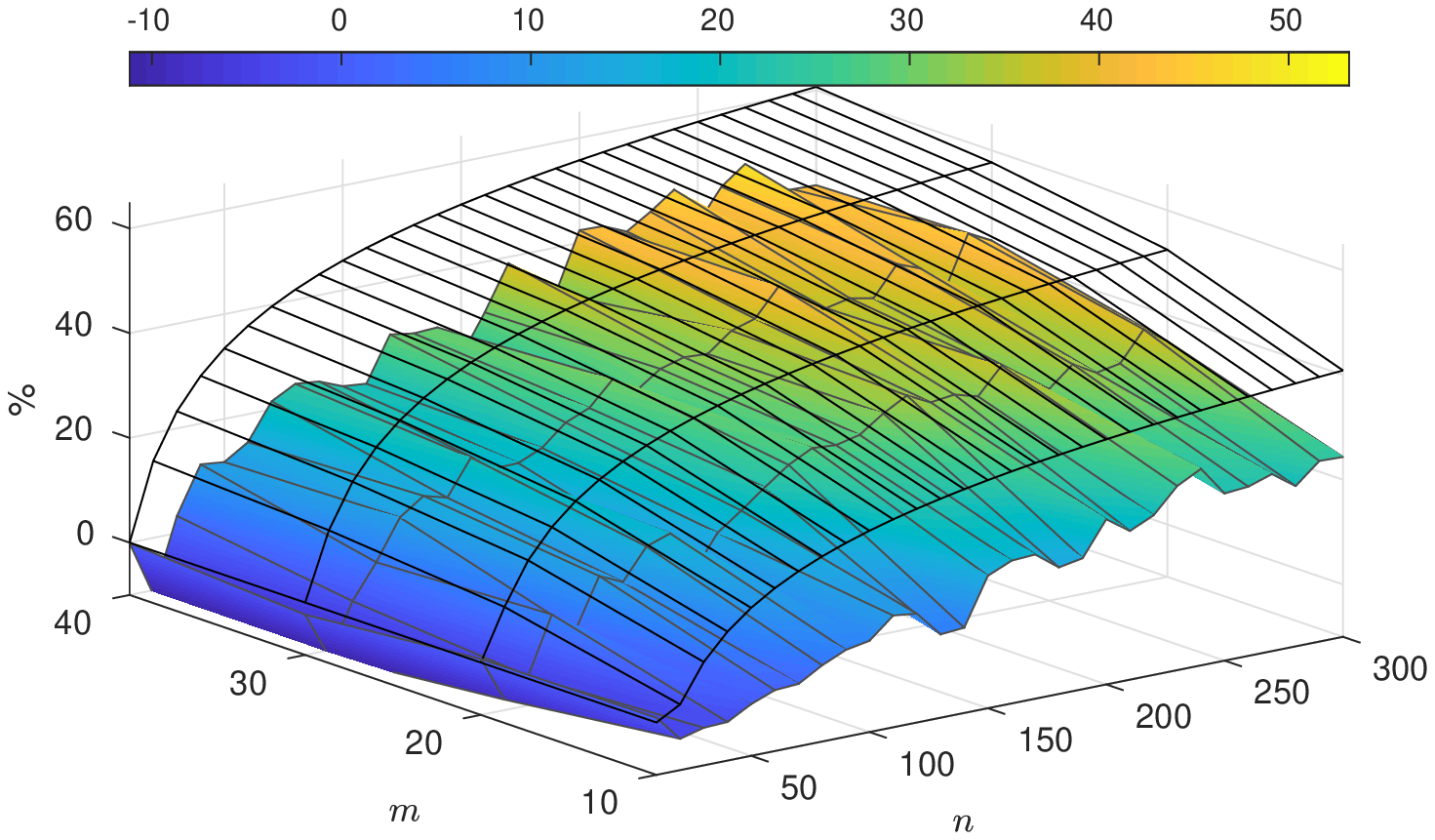}
  \end{subfigure}
  \hspace*{0.01\linewidth}
  \begin{subfigure}[b]{0.485\linewidth}
    \centering
    \includegraphics[width=\linewidth]{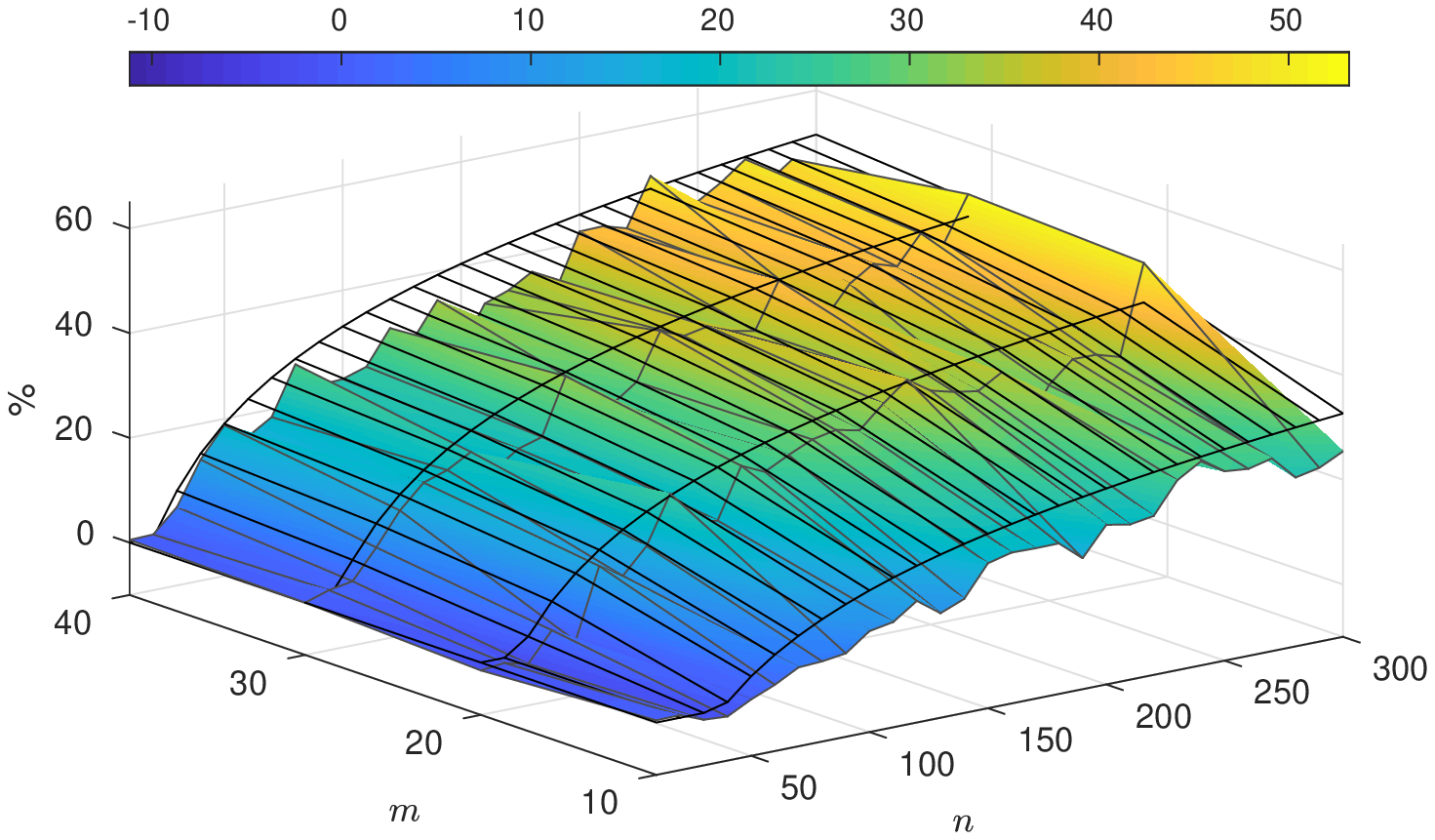}
  \end{subfigure}
    \caption{For the metrics $\nu_{\rm rows}$ (left) and $\nu_{\rm flops}$ (right), we compare the additional cost $N\big(\{\set{J}\},\epsilon,\beta\big)/N\big(\setbf{P}^\star,\greekbf{\epsilon}^\star, \tfrac{\beta}{|\setbf{P}^\star|}\mathds{1}\big) -1$ in \% (grid) to the additional actual solver time $\big(\text{solver time } \PSP{\{\set{J}\}}\big)/\big(\text{solver time } \PSP{\setbf{P}^\star}\big) -1$ in \% (surface), taking the median over $20$ random instances.
    }
    \label{ccp:example:fig_ratios}
\end{figure}

When the cost metric $\nu$ is the sample complexity, Algorithm~1 returns the trivial partition $\{\set{J}\}$ for all instances. However, when $\nu$ is selected to be the number of constraint rows, i.e., $\nu_{\rm rows}(\set{A}) := |\set{A}|$, the number of non-zero elements $\nu_{\rm nnzs}$, or the number of FLOPs $\nu_{\rm flops}$, then Algorithm~1 returns non-trivial partitions for most instances. In these cases, using $\setbf{P}^\star$, $\greekbf{\epsilon}^\star$ leads to a reduction of the cost compared to using the trivial partition. Moreover, the ratio of this reduction grows with growing number of products $m$ or machines $n$.
In \reffig{ccp:example:fig_ratios}, we consider the two cost metrics $\nu_{\rm rows}$ and $\nu_{\rm flops}$, and compare $\PSP{\{\set{J}\}}$ to $\PSP{\setbf{P}^\star}$. We also report the additional computational effort needed for solving $\PSP{\{\set{J}\}}$ compared to $\PSP{\setbf{P}^\star}$ based on median solver times. For $\nu_{\rm flops}$, and for $n=300$ and $m=40$ the median solver times, of $20$ instances, are $8$~seconds for $\PSP{\{\set{J}\}}$ and $5.6$~seconds for $\PSP{\setbf{P}^\star}$.
The metric $\nu_{\rm flops}$ is able to accurately predict the improvement in computational effort and outperforms $\nu_{\rm rows}$ in this respect. This illustrates that the cost metric $\nu_{\rm flops}$ is indicative of the actual computational effort.
Moreover, when $\nu_{\rm flops}$ is used to compute $\setbf{P}^\star$, $\greekbf{\epsilon}^\star$, the computational effort of $\PSP{\setbf{P}^\star}$ is reduced compared to $\nu_{\rm rows}$ in 80\% of instances. In contrast, the difference between using $\nu_{\rm nnzs}$ or $\nu_{\rm flops}$ is negligible.
Note that while \refalg{ccp:algorithm_paritioning} uses the explicit bound \eqref{ccp:eq:explicit_lb_extended}, we used the implicit sample bound from \eqref{ccp:eqn:implicit} to construct $\PSP{\{\set{J}\}}$ and $\PSP{\setbf{P}^\star}$.
Despite this, the solver times as seen in \reffig{ccp:example:fig_ratios} are accurately predicted by the cost metric $\nu_{\rm flops}$ used in \refalg{ccp:algorithm_paritioning}.
Finally, \reffig{ccp:example:fig_ratios} illustrates that the proposed partitioning scheme can be used to substantially reduce the computational effort required for the scenario approach.

\subsection{Multi-agent formation planning}\label{ccp:case_study:agents}
We consider the task of planning the motion of $n$ agents to a target formation consisting of $n$ formation points $p^j \in \reals{2}$ for $j=1,\ldots,n$. Agent $i$ is modeled as a two-dimensional double integrator discretized at $0.1$~seconds: 
\begin{equation*}
 \begin{bmatrix}\dot{x}^{i} \\ \dot{y}^{i} \end{bmatrix} = \begin{bmatrix}v^{i}_x \\ v^{i}_y\end{bmatrix}, \text{ and } \begin{bmatrix}\dot{v}^{i}_x \\ \dot{v}^{i}_y\end{bmatrix} = u^{i} + w^{i}\,,
\end{equation*}
where $(x^{i}_k, y^{i}_k, v^{i}_{x,k}, v^{i}_{y,k}) \in \reals{4}$, are horizontal and vertical positions, and forward and lateral velocities of agent $i$ at time $k$, respectively.
The control inputs $u^{i}_k \in \reals{n_u}$, with $n_u := 2$, are the forward and lateral accelerations applied to agent $i$ from time $k$ to $k+1$, and are affected by an additive disturbance $w^{i}_k \in \reals{2}$ with $w^{i}_{kj} \sim \uniform [-0.25, 0.25]$, modeling actuator uncertainty or adversarial control action.
The planning horizon is $N=20$.

\begin{figure*}[ht]
 \centering
 \begin{minipage}{0.38\linewidth}
   \centering%
   \begin{subfigure}[b]{0.485\linewidth}
     \includetikzfigure{trajectories_example2_3_3flops}
     \caption{$n=3$ agents}
   \end{subfigure}
   \hspace*{0.01\linewidth}
   \begin{subfigure}[b]{0.47\linewidth}
     \includetikzfigure{trajectories_example2_5_5flops}
     \caption{$n=5$ agents}
   \end{subfigure}
   \caption{100 Random realizations } \label{ccp:example2:task}
 \end{minipage}
 \hspace*{0.01\linewidth}
 \begin{minipage}{0.56\linewidth}
   \centering%
   \begin{subfigure}[b]{0.55\linewidth}
     \centering
     \includetikzfigure{empiricalViolations_example2_5_5flops}
     \caption{Empirical violations of $f_{\set{P}_1^\star}$,$f_{\set{P}_2^\star}$, and $f$ for \num{10000} new realizations.}
     \label{ccp:example2:empirical_violations}
   \end{subfigure}
   \hspace*{0.025\linewidth}
   \begin{subfigure}[b]{0.40\linewidth}
     \centering%
     \includetikzfigure{solvertimes_example2_5_5flops}
     \caption{Solver times for $\PSP{\{\set{J}\}}$ and $\PSP{\setbf{P}^\star}$.}
     \label{ccp:example2:solvertime}
   \end{subfigure}
   \caption{The median ({\color{red}\textbf{red}}), and the 25$^{th}$ and 75$^{th}$ percentiles ({\color{blue}\textbf{blue}}).}
 \end{minipage}
\end{figure*}

The formation points are the vertices of a regular $n$-gon with circumradius $0.5$ centered at the average final position of the agents $z := \frac{1}{n}\sum_{i=1}^n (x^{i}_N, y^{i}_N)$.
For agent $i$ assigned to the $j$-th vertex $p^j$ we require that $\|(x^{i}_N, y^{i}_N) - (z+p^{j}) \|_\infty \leq 0.35$.
We furthermore have input bounds $-\mathds{1} \leq \tfrac{1}{2}(u^{i}_k + w^{i}_k) \leq \mathds{1}$ for $i=1,\ldots,n$ and $k=0,\ldots,N-1$.
For each agent $i$ the objective is to minimize the infinity-norm distance $\|(x^{i}_N, y^{i}_N) - (z+p^{j}) \|_\infty$ to its assigned vertex $p^j$ and a quadratic cost on the control inputs.
We define $\mathbf{u}^{i} := (u^{i}_0, \ldots, u^{i}_{N-1}) \in \reals{n_uN}$ and $\mathbf{u} := (\mathbf{u}^{1},\ldots,\mathbf{u}^{n})$, and analogously $\mathbf{w}^{i}$, $\mathbf{w}$. Note that $z$ and $(x_N^{i},y_N^{i})$ can be written in terms of $\mathbf{u}$, $\mathbf{w}$ and the known initial state $\theta^{i} \in \reals{4}$. This gives rise to the following chance constrained program:
\begin{align}\label{ccp:example2:robust_program}
 \min_{\mathbf{u}, \greekbf{\tau}, \mathbf{t}, \mathbf{y}} \:\:&  \frac{1}{4nN}\|\mathbf{u}\|_2^2 + \sum_{i=1}^n \mathds{1}^\transp (\greekbf{\tau}^{i} - \mathbf{t}^{i}) \\
 \suchthat \:& \mathbb{P}_{\mathbf{w}} \big[ \underbrace{\begin{bmatrix}
   g_1(\mathbf{u},\greekbf{\tau},\mathbf{w}) \\
   g_2(\mathbf{u},\mathbf{w})
 \end{bmatrix}}_{=: f(\mathbf{u},\greekbf{\tau},\mathbf{w})} \leq 0 \big] \geq (1-\epsilon) \,, \nonumber\\
 & \greekbf{\tau}^{i} \leq 0.35\mathds{1} + M(\mathds{1}-\mathbf{y}^{i}) && \hspace*{-2em} i=1,\ldots,n \nonumber\,,\\
 & \mathbf{t}^{i} \leq M(\mathds{1}-\mathbf{y}^{i}) && \hspace*{-2em} i=1,\ldots,n\,, \nonumber\\
 & 0 \leq \mathbf{t}^{i} \leq \greekbf{\tau}^{i} && \hspace*{-2em} i=1,\ldots,n\,, \nonumber\\
 & \sum_{i=1}^n\nolimits \mathbf{y}^{i} \leq \mathds{1}\,,\: \mathds{1}^\transp \mathbf{y}^{i} = 1 && \hspace*{-2em} i=1,\ldots,n\,, \nonumber\\
 & \mathbf{y}^{i} \in \binaries{m} && \hspace*{-2em} i=1,\ldots,n\,, \nonumber
\end{align}
where $g_{1} := (g_{1}^{1,1}, \ldots, g_{1}^{n,n})$ encodes the formation constraint, with $g_{1}^{i,j}(\mathbf{u},\greekbf{\tau},\mathbf{w}) := \begin{bsmallmatrix} x^{i}_N - (z+p^j) - \tau^{i}_j\mathds{1} \\ (z+p^j) - x^{i}_N - \tau^{i}_j\mathds{1} \end{bsmallmatrix} \in \reals{4}$, and $g_{2}(\mathbf{u},\mathbf{w}) := \begin{bsmallmatrix}\mathbf{u}-\mathbf{w} - 2\mathds{1} \\ \mathbf{w}-\mathbf{u} -2\mathds{1}\end{bsmallmatrix} \in \reals{2nn_uN}$ encodes the input constraints.
The auxiliary variable $\greekbf{\tau} := (\greekbf{\tau}^{1},\ldots, \greekbf{\tau}^{n})$ with $\greekbf{\tau}^{i} \in \reals{n}$ is bounded by the infinity-norm distance of agent $i$ to each of the vertices of the $n$-gon.
Each of the binary vectors $\mathbf{y}^{i} \in \binaries{n}$ has exactly one non-zero entry whose index $j$ indicates that agent $i$ is assigned to vertex $p^j$.
The $j$-th entry $t^{i}_j$ of the additional auxiliary variable $\mathbf{t} := (\mathbf{t}^{1},\ldots, \mathbf{t}^{n})$ with $\mathbf{t}^{i} \in \reals{n}$ is zero if $y^{i}_j = 1$ and equal to $\tau^{i}_j$ otherwise. Therefore $\mathds{1}^\transp (\greekbf{\tau}^{i} - \mathbf{t}^{i})$ encodes the infinity-norm cost of the deviation of agent $i$ from its assigned vertex $p^j$. This is encoded via linear constraints on $\greekbf{\tau}$, $\mathbf{t}$ and $\mathbf{y}$ using the Big-M formulation \cite{bemporad1999} with a large enough $M \in \posreals{}$.

Given $\epsilon = 0.1$ and $\beta = 10^{-3}$, we look for approximate solutions to \eqref{ccp:example2:robust_program} via the scenario approach.
We use the computational cost metric $\nu_{\rm flops}$.
For $n=5$ agents, \refalg{ccp:algorithm_paritioning} returns the partition $\setbf{P}^\star = \{\set{P}_1^\star, \set{P}_2^\star\}$ and $\greekbf{\epsilon}^\star = (0.0313, 0.0687)$ with $\set{P}_1^\star = \{1,\ldots,4n^2\}$ and $\set{P}_2^\star = \{4n^2+1,\ldots,4n^2+2nn_uN\}$. Note that the number of possible partitions is $B_r = B_{500} \approx 10^{855}$, where $B_r$ is the $r$-th Bell number \cite{bell1934}.
The constraint functions of the scenario programs $\PSP{\{\set{J}\}}$ and $\PSP{\setbf{P}^\star}$ require a total number of \unit[30]{MFLOPs} and \unit[11]{MFLOPs}, respectively.
For $n=3$ and $n=5$, the solution of $\PSP{\setbf{P}^\star}$ is illustrated in \reffig{ccp:example2:task} for $100$ new random realizations.

To evaluate the theoretical guarantees for $n=5$, we solved $\PSP{\setbf{P}^\star}$ for $200$ instances of the multisample.
Computations were carried out on an Intel~i7 CPU with $4$ cores at $\unit[2.8]{GHz}$ and $\unit[8]{GB}$ of memory, using YALMIP \cite{lofberg2004} and CPLEX \cite{cplex} as optimization problem solver.
For each solution, we evaluated the constraint $f$ over \num{10000} new realizations of the uncertain parameters $\mathbf{w}$. \reffig{ccp:example2:empirical_violations} shows the distribution of the obtained empirical violations $\hat{\epsilon}_1$, $\hat{\epsilon}_2$, and $\hat{\epsilon}$ of $f_{\set{P}^\star_1}$, $f_{\set{P}^\star_2}$, and $f$, respectively. As expected by \reflem{ccp:lemma:sampling_lemma} and \refthm{ccp:sec2:thm1}, $\hat{\epsilon}_1, \hat{\epsilon}_2$ and $\hat{\epsilon}$ satisfy the parameters $\epsilon_1^\star, \epsilon_2^\star$ and $\epsilon$, i.e., the desired probabilistic guarantees are satisfied.
In \reffig{ccp:example2:solvertime} we further show that $\setbf{P}^\star$ leads to reduced solver times compared to the trivial partition $\{\set{J}\}$.

\section{Conclusions}\label{ccp:sec:conclusions}
We have proposed an approach to reduce the computational effort required for solving scenario programs.
We have formulated a nonlinear optimization problem to find the partition that minimizes the computational cost of the scenario program. We proved that the support rank is a monotone submodular function of the constraint row indices and have proposed a polynomial-time algorithm which finds suboptimal solutions to this partitioning problem. For the submodular case, we provided approximation guarantees for the proposed algorithm.
Finally, we have demonstrated that our approach leads to computational advantages in two case studies from production and multi-agent planning. Encouraged by the obtained results, extending the approximation guarantees of \textsc{SubMP} to $\mu$-weakly submodular functions is subject of our current research.

\bibliographystyle{plain}        % Include this if you use bibtex 
\bibliography{mybibfile}           % and a bib file to produce the 
                                 % correct style is generated by
                                
\appendix

\section{Auxiliary proofs and results}\label{ccp:apdx:proofs}
\vspace{-\parskip}
\def\Elproofname{\protect{PROOF of \reflem{ccp:lem:subspace}.}}
\begin{pf} \label{ccp:pf:subspace}
  For a function $g : \reals{n} \times \binaries{b} \times \setbf{W} \rightarrow \reals{}$, $\set{L}_g$ is defined in \cite[Eqn.~(3.7)]{schildbach2013} and is the set of linear subspaces of $\reals{n}$ not constrained by $g$, i.e.,
  \begin{equation}
    \set{L}_g := \hspace*{-1em}\bigcap_{(\mathbf{x},\mathbf{y},\mathbf{w}) \in \setbf{Z}} \hspace*{-1em} \big\{L \subseteq \set{L} \sep{\big} L \subseteq F_g(\mathbf{x},\mathbf{y},\mathbf{w}) \big\}\,, \label{eqn:Lg}
  \end{equation}
  where $\setbf{Z} := \setbf{V} \times \setbf{W}$ and $F_g(\mathbf{x},\mathbf{y},\mathbf{w}) := \{\greekbf{\xi} \in \reals{n} \sep{} g(\mathbf{x}+\greekbf{\xi},\mathbf{y},\mathbf{w}) = g(\mathbf{x},\mathbf{y},\mathbf{w})\}$, and $\set{L}$ is the set of all linear subspaces of $\reals{n}$.
  \vspace*{-\parskip}
  \begin{enumerate}[label=\roman*),wide,nosep]
    \item Consider $g,h : \reals{n} \times \binaries{b} \times \setbf{W} \rightarrow \reals{}$ with $g(\mathbf{x},\mathbf{y},\mathbf{w}) \geq h(\mathbf{x},\mathbf{y},\mathbf{w})$ for all $(\mathbf{x},\mathbf{y},\mathbf{w}) \in \setbf{Z}$, with $g,h$ proper and $h$ convex in $\mathbf{x}$. 
    We will show that $\set{L}_g \subseteq \set{L}_h$. This then gives the result using $h(\mathbf{x},\mathbf{y},\mathbf{w}) := \max_{j\in\set{A}} f_j(\mathbf{x},\mathbf{y},\mathbf{w})$ and $g(\mathbf{x},\mathbf{y},\mathbf{w}) := \max_{j\in\set{B}} f_j(\mathbf{x},\mathbf{y},\mathbf{w})$ since $\set{A} \subseteq \set{B}$. By definition of $f$, both $g$ and $h$ are proper and convex in $\mathbf{x}$.
    We consider any non-zero dimensional linear subspace $L \in \set{L}_g$ and any line $\{\lambda \mathbf{d} \sep{} \lambda \in \reals{} \} \subseteq L$, with $\mathbf{d} \in \reals{n}\setminus\{0\}$. By definition of $\set{L}_g$, see \eqref{eqn:Lg}, $g(\mathbf{x}+\lambda \mathbf{d},\mathbf{y},\mathbf{w}) = g(\mathbf{x},\mathbf{y},\mathbf{w})$ for all $(\mathbf{x},\mathbf{y},\mathbf{w}) \in \setbf{Z}$ and all $\lambda \in \reals{}$.
    This implies that for any fixed $(\bar{\mathbf{x}},\bar{\mathbf{y}},\bar{\mathbf{w}}) \in \setbf{Z}$ we have that $h(\bar{\mathbf{x}}+ \lambda \mathbf{d},\bar{\mathbf{y}},\bar{\mathbf{w}}) \leq g(\bar{\mathbf{x}}+ \lambda \mathbf{d},\bar{\mathbf{y}},\bar{\mathbf{w}}) = g(\bar{\mathbf{x}},\bar{\mathbf{y}},\bar{\mathbf{w}})$, for all $\lambda \in \reals{}$.
    Due to convexity of $h$ in $\lambda$ (and $g,h$ proper) this is only true if $h(\bar{\mathbf{x}},\bar{\mathbf{y}},\bar{\mathbf{w}}) = h(\bar{\mathbf{x}}+ \lambda \mathbf{d},\bar{\mathbf{y}},\bar{\mathbf{w}})$ for all $\lambda \in \reals{}$.
    Otherwise, convexity of $h$ would imply the existence of a $\mu \in \reals{}$ such that $h(\bar{\mathbf{x}}+ \mu \mathbf{d},\bar{\mathbf{y}},\bar{\mathbf{w}})$ is arbitrarily large, i.e., $h(\bar{\mathbf{x}}+ \mu \mathbf{d},\bar{\mathbf{y}},\bar{\mathbf{w}}) > g(\bar{\mathbf{x}}+ \mu \mathbf{d},\bar{\mathbf{y}},\bar{\mathbf{w}}) = g(\bar{\mathbf{x}},\bar{\mathbf{y}},\bar{\mathbf{w}})$, which leads to a contradiction.
    Since $(\bar{\mathbf{x}},\bar{\mathbf{y}},\bar{\mathbf{w}})$ was chosen arbitrarily, and the choice of $\mathbf{d}$ and $L$ was also arbitrary, it follows that $L \in \set{L}_h$ and thus $\set{L}_g \subseteq \set{L}_h$, which gives the result. \hfill\(\qed\)
    \item Consider $g_1,g_2,h : \reals{n} \times \binaries{b} \times \setbf{W} \rightarrow \reals{}$ with $g_1,g_2,h$ proper and $g_1,g_2$ convex in their first argument. We define $h(\mathbf{x},\mathbf{y},\mathbf{w}) := \max\{g_1(\mathbf{x},\mathbf{y},\mathbf{w}), g_2(\mathbf{x},\mathbf{y},\mathbf{w})\}$. 
    We will show that $\set{L}_h = \set{L}_{g_1} \cap \set{L}_{g_2}$, which gives the result using $g_1(\mathbf{x},\mathbf{y},\mathbf{w}) := \max_{\hat{\jmath}\in\set{A} \cup \{j\}} f_{\hat{\jmath}}(\mathbf{x},\mathbf{y},\mathbf{w})$, $g_2(\mathbf{x},\mathbf{y},\mathbf{w}) := \max_{\hat{\jmath}\in\set{B}} f_{\hat{\jmath}}(\mathbf{x},\mathbf{y},\mathbf{w})$, which using $\set{A} \subseteq \set{B}$ implies that $h(\mathbf{x},\mathbf{y},\mathbf{w}) = \max_{\hat{\jmath}\in\set{B} \cup \{j\}} f_{\hat{\jmath}}(\mathbf{x},\mathbf{y},\mathbf{w})$. The functions $g_1$, $g_2$, and $h$ are proper and convex in their first argument by definition of $f$.
    Consider any linear subspace $L \in \set{L}_{g_1} \cap \set{L}_{g_2}$ and line $\{\lambda \mathbf{d} \sep{} \lambda \in \reals{} \} \subseteq L$, with $\mathbf{d} \in \reals{n}\setminus\{0\}$. By definition of $\set{L}_{g_1}$ and $\set{L}_{g_2}$, see \cite[Eqn.~(3.6)]{schildbach2013}, we have $g_1(\mathbf{x}+\lambda \mathbf{d},\mathbf{y},\mathbf{w}) = g_1(\mathbf{x},\mathbf{y},\mathbf{w})$ and $g_2(\mathbf{z}+\lambda \mathbf{d},\mathbf{y},\mathbf{w}) = g_2(\mathbf{x},\mathbf{y},\mathbf{w})$ for all $(\mathbf{x},\mathbf{y},\mathbf{w}) \in \setbf{Z}$ and all $\lambda \in \reals{}$. It follows that for any $(\mathbf{x},\mathbf{y},\mathbf{w}) \in \setbf{Z}$ and any $\lambda \in \reals{}$, we have that $h(\mathbf{x}+\lambda \mathbf{d},\mathbf{y},\mathbf{w}) = \max\{g_1(\mathbf{x}+\lambda \mathbf{d},\mathbf{y},\mathbf{w}), g_2(\mathbf{x}+\lambda \mathbf{d},\mathbf{y},\mathbf{w})\} = \max\{g_1(\mathbf{x},\mathbf{y},\mathbf{w}), g_2(\mathbf{x},\mathbf{y},\mathbf{w})\} = h(\mathbf{x},\mathbf{y},\mathbf{w})$. This implies that $L \in \set{L}_h$ for any $L$ and gives $\set{L}_{g_1} \cap \set{L}_{g_2} \subseteq \set{L}_h$. Applying \ref{ccp:lem:subspace:mon} with convexity of $g_1,g_2$ we also get $\set{L}_h \subseteq \set{L}_{g_1} \cap \set{L}_{g_2}$, which implies $\set{L}_h = \set{L}_{g_1} \cap \set{L}_{g_2}$, giving the result. \hfill\(\qed\)
  \end{enumerate}
\end{pf}
   \def\Elproofname{\protect{PROOF.}}

We will use that a set function $g: 2^\mathcal{J} \rightarrow \reals{}$ is modular if and only if $g(\set{A}) = \sum_{j \in \set{A}} \big(g(\{j\}) - \nu(\emptyset) \big)+ g(\emptyset) $ for all $\set{A} \subseteq \set{J}$.
\begin{proposition}
\label{prop:composition_with_concave}
Let $\gamma : 2^\mathcal{J} \rightarrow \posreals{}$ be a monotone, $\mu$-weakly submodular set function, and $g: \mathbb{R}\rightarrow \mathbb{R}$ be a non-decreasing \emph{concave} function. Then, the composition $f= g \circ \gamma$ is $\mu$-weakly submodular.
\end{proposition}
\vspace*{-\parskip}
\begin{pf}
\sloppy
Concavity and non-decreasingness of $g$ imply that $\forall x,y,z,w \in \reals{}$ with $x < y, z<w, x \leq z$, we have $\big(g(y)- g(x)\big)/\big(g(w) - g(z)\big) \geq (y-x)/(w-z)$. Therefore, for any pair $\set{A},\set{B} \subseteq \mathcal{J}$ with $\set{A} \cap \set{B} = \emptyset$ and $\set{B} \neq \emptyset$, 
\begin{align*}
\frac{\sum_{j\in \set{B}} f(j \sep{} \set{A})}{ f(\set{B}\sep{}\set{A})} & = \frac{\sum_{j\in \set{B}} g(\gamma(\set{A} \cup j)) -  g(\gamma(\set{A}))}{ g(\gamma(\set{A} \cup \set{B})) -  g(\gamma(\set{A}))}\\ 
& \geq \frac{\sum_{j\in \set{B}} \gamma(\set{A} \cup j) -  \gamma(\set{A})}{ \gamma(\set{A} \cup \set{B}) -  \gamma(\set{A})} \geq \mu \, .
\end{align*}
The first inequality follows since $g$ is concave and non-decreasing and since for any $j\in \set{B}$ with $\gamma(\set{A}\cup j) \neq \gamma(\set{A})$, we have $\gamma(\set{A} \cup \set{B}) \geq \gamma(\set{A} \cup j)  > \gamma(\set{A})$. The second inequality is due to $\gamma$ being $\mu$-weakly submodular.  
For $\set{B}=\emptyset$, $\mu$-weak submodularity is trivially satisfied.
   \hfill\(\qed\)
\end{pf}

\begin{lemma}
\label{lemma:product_modular_weak_submodular}
Let $f,g : 2^\mathcal{J} \rightarrow \posreals{}$ be two monotone set functions, with $f$ \emph{submodular}, and $g$ \emph{modular} such that $g(\set{A}) > 0$ for all $\set{A} \subseteq \set{J}$, with $\set{A} \neq \emptyset$. Then, there exists a $\mu > 0$ such that the product $\gamma(\set{A}) := f(\set{A}) g(\set{A})$ is $\mu$-weakly submodular.
\end{lemma}
\vspace*{-\parskip}
\begin{pf}
Note that by definition of $\gamma$ and modularity of $g$, for any $\set{A},\set{B}\subseteq \mathcal{J}$ we have that $\gamma(\set{B}\sep{}\set{A}) =  f(\set{A}\cup \set{B})g(\set{A}\cup \set{B})-  f(\set{A})g(\set{A}) = f(\set{B}\sep{}\set{A})g(\set{A}\cup \set{B}) + f(\set{A})\big(g(\set{B}) - g(\emptyset)\big)$. Consider now any pair $\set{A},\set{B} \subseteq \mathcal{J}$ with $\set{A} \cap \set{B} = \emptyset$. 
Note that for $\set{B} = \emptyset$, $\mu$-weak submodularity trivially holds for any $\mu$.
For $\set{B} \neq \emptyset$ we have: 
\begin{align*}
&\sum_{j\in \set{B}} \gamma(j|\set{A}) = \hspace*{-0.2em}\sum_{j\in \set{B}} \hspace*{-0.2em} \Big( f(j|\set{A}) g(\set{A}\cup j) + f(\set{A})\big(g(j) - g(\emptyset)\big) \Big) \\
& \geq \min_{j \in \set{B}}g(\set{A}\cup j)\sum_{j\in \set{B}} f(j\sep{}\set{A})  + f(\set{A})\big(g(\set{B}) - g(\emptyset)\big)  \\
& \geq  \min_{j \in \set{B}}g(\set{A}\cup j)f(\set{B}\sep{}\set{A})  + f(\set{A})\big(g(\set{B}) - g(\emptyset)\big)  \\
& \geq \underbrace{\frac{\min_{j \in \set{B}}g(\set{A}\cup j)}{g(\set{A}\cup \set{B})}}_{:= \mu(\set{A},\set{B})} \Big( f(\set{B}\sep{}\set{A}) g(\set{A}\cup \set{B})\\[-2.5em]
& \phantom{\geq \frac{\min_{j \in \set{B}}g(\set{A}\cup j)}{g(\set{A}\cup \set{B})} \Big(} \: +  f(\set{A})\big(g(\set{B}) - g(\emptyset)\big) \Big)\\
&  =\mu(\set{A},\set{B}) \gamma(\set{B}\sep{}\set{A}) \,.
\end{align*}
The first inequality follows by modularity of $g$ and non-negativity of $f$ and $g$, and the second inequality follows from ($1$-weak) submodularity of $f$. The last inequality is due to $0 < \mu(\set{A},\set{B}) \leq 1$, which holds because of $g(\set{A}) > 0$ for all $\set{A} \neq \emptyset$ and monotonicity of $g$. From monotonicity of $g$, it then follows that $\gamma$ is  $\mu$-weakly submodular with 
\begin{equation*}
\mu := \min_{\substack{\set{A}, \set{B} \subseteq \mathcal{J} \\ \set{A} \cap \set{B}= \emptyset }}{  \mu(\set{A},\set{B}) } = \mu(\emptyset,\set{J}) > 0 \,. \qed
\end{equation*}
\end{pf}

\section{Extension to mixed-integer convex case} \label{ccp:apdx:mixed-integer}
For notational simplicity, we define $V^\star[\mathbb{K}] := V(\mathbf{x}^\star[\mathbb{K}],\mathbf{y}^\star[\mathbb{K}])$ and $V_i^\star[\mathbb{K}] := V_i(\mathbf{x}^\star[\mathbb{K}],\mathbf{y}^\star[\mathbb{K}])$.
\begin{lemma}\label{ccp:lemma:sampling_lemma}
Consider $\MCP$. Let $(\mathbf{x}^\star[\mathbb{K}], \mathbf{y}^\star[\mathbb{K}])$ be the unique optimizer of the scenario program $\SP$. Then, for $i=1,\ldots, P$ and parameters $\epsilon_i \in (0,1)$, 
\begin{equation*}
\mathbb{P}^K [V_i^\star[\mathbb{K}] > \epsilon_i ] \\ \leq 2^b \sum_{l=0}^{\rho(\set{P}_i)-1} \hspace{-0.4em} \binom{K_i}{l} \epsilon_i^l (1-\epsilon_i)^{K_i-l} \,.
\end{equation*}
\end{lemma}
\vspace*{-\parskip}
\begin{pf}
The proof follows similarly to \cite{schildbach2013}. Given any $\hat{\imath} \in \set{I}$, we define $\set{I}^{\minus\hat{\imath}} := \set{I}\setminus \{\hat{\imath}\}$. First, we  will upper bound the probability of $V_{\hat{\imath}}^\star[\mathbb{K}]$ being larger than $\epsilon_{\hat{\imath}}$ conditioned on the samples $\mathbf{W}^{\minus\hat{\imath}} := \{ \mathbf{w}^{(i,j)} \sep{} i\in \set{I}^{\minus\hat{\imath}} , j \in \set{K}_i \}$ drawn for the other constraints, i.e.,
\begin{equation}\label{ccp:eq:condit_prob}
\mathbb{P}^K \big[ V_{\hat{\imath}}^\star[\mathbb{K}] > \epsilon_{\hat{\imath}} \sep{\big} \mathbf{W}^{\minus\hat{\imath}} \big]\,.
\end{equation}
Then we marginalize out $\mathbf{W}^{\minus\hat{\imath}}$ to obtain the result.

To upper bound \eqref{ccp:eq:condit_prob}, we reason about the probability
\begin{equation} \label{ccp:eq:thm:sampling:proofprob}
  \mathbb{P}^{K_{\hat{\imath}}} \big[ V_{\hat{\imath}}^\star[\mathbb{K}] > \epsilon_{\hat{\imath}} \sep{\big} \overbar{\mathbf{w}}^{(i,j)} , \forall\: i \in \set{I}^{\minus\hat{\imath}}, j \in \set{K}_i \big]\,,
\end{equation} 
where $\{ \overbar{\mathbf{w}}^{(i,j)} \sep{} i\in \set{I}^{\minus\hat{\imath}} , j \in \set{K}_i \}$ is a realization of the samples $\mathbf{W}^{\minus\hat{\imath}}$, and $\mathbb{P}^{K_{\hat{\imath}}}$ is the probability distribution of the samples $\mathbf{W}^{\hat{\imath}} := \{ \mathbf{w}^{(\hat{\imath},j)} \sep{} j \in \set{K}_{\hat{\imath}} \}$.
We first show that
\begin{equation} \label{ccp:sec:2:eq:conditional_prob_bound}
  \eqref{ccp:eq:thm:sampling:proofprob} \leq 2^b\sum_{l=0}^{\rho(\set{P}_{\hat{\imath}})-1} \binom{K_{\hat{\imath}}}{l} \epsilon_{\hat{\imath}}^l (1-\epsilon_{\hat{\imath}})^{K_{\hat{\imath}}-l}\,.
\end{equation}
For fixed realizations $\overbar{\mathbf{w}}^{(i,j)}$, $\SP$ becomes: 
\begin{align} \label{ccp:sec:2:eq:first_class_conditional}
  \min_{(\mathbf{x},\mathbf{y}) \in \setbf{V}} \:&  c^\transp \mathbf{x}  \\[-0.4em]
  \suchthat \:& f_{\hat{\imath}}(\mathbf{x}, \mathbf{y},\mathbf{w}^{(\hat{\imath},j)}) \leq 0 \quad \forall j  \in \set{K}_{\hat{\imath}}\,,\nonumber \\
  \:& f_{\minus\hat{\imath}}(\mathbf{x},\mathbf{y}) := \hspace*{-0.5em} \max_{i \in \set{I}^{\minus\hat{\imath}},\,j\in\set{K}_i} \hspace*{-0.5em}  f_i(\mathbf{x}, \mathbf{y},\overbar{\mathbf{w}}^{(i,j)} )  \leq 0 \,, \nonumber
 \end{align}
where $f_{\minus\hat{\imath}}$ is a deterministic mixed-integer constraint that is convex in $\mathbf{x}$ for each $\mathbf{y}$, as the point-wise maximum of convex functions of $\mathbf{x}$ for each $\mathbf{y}$. \eqref{ccp:sec:2:eq:first_class_conditional} is the scenario version of a CCP, where all except the $\hat{\imath}$-th constraint row have been made deterministic by substituting the realizations $\overbar{\mathbf{w}}^{(i,j)}$. Moreover, \eqref{ccp:sec:2:eq:first_class_conditional} is convex for each configuration of binary variables $\mathbf{y}$. Therefore, it falls into the class of problems considered in \cite{esfahani2015}.
Let $(\mathbf{x}^\star[K_{\hat{\imath}}],\mathbf{y}^\star[K_{\hat{\imath}}])$ be the optimizer of \eqref{ccp:sec:2:eq:first_class_conditional}. Note that $(\mathbf{x}^\star[K_{\hat{\imath}}],\mathbf{y}^\star[K_{\hat{\imath}}])$ is random and its distribution is a function of $\mathbb{P}^{K_{\hat{\imath}}}$ depending only on the samples in $\mathbf{W}^{\hat{\imath}}$. Following the probabilistic arguments in \cite[Thm.~4.1]{esfahani2015}, we obtain \eqref{ccp:sec:2:eq:conditional_prob_bound}. Since the preceding discussion holds for any set of realizations $\overbar{\mathbf{w}}^{(i,j)}$ of $\mathbf{W}^{\minus\hat{\imath}}$, it immediately follows that
\begin{align} \label{ccp:sec:2:eq:conditional_prob_bound_K}
\eqref{ccp:eq:condit_prob} \leq 2^b\sum_{l=0}^{\rho(\set{P}_{\hat{\imath}})-1} \binom{K_{\hat{\imath}}}{l} \epsilon_{\hat{\imath}}^l (1-\epsilon_{\hat{\imath}})^{K_{\hat{\imath}}-l} \,.
\end{align}
We obtain the result of \reflem{ccp:lemma:sampling_lemma} by marginalizing out the conditioned samples. The right hand side of \eqref{ccp:sec:2:eq:conditional_prob_bound_K} is independent from 
$\mathbf{W}^{\minus\hat{\imath}}$. Therefore, with $ \mathbf{W}^{\minus\hat{\imath}}$ distributed according to $\mathbb{P}^{K_{\minus\hat{\imath}}} :=  \mathbb{P}^{K_1 + \dots K_{\hat{\imath}-1} + K_{\hat{\imath}+1} + \dots K_P}$ we have
\begin{align*}
  \mathbb{P}^K \big[ V_{\hat{\imath}}^\star[\mathbb{K}] > \epsilon_{\hat{\imath}} \big] &= \int_{\mathbf{W}^{\minus\hat{\imath}}} \mathbb{P}^K \big[ V_{\hat{\imath}}^\star[\mathbb{K}] > \epsilon_{\hat{\imath}} \sep{\big} \mathbf{W}^{\minus\hat{\imath}} \big] \mathbb{P}^{K_{\minus\hat{\imath}}}[ d\, \mathbf{W}^{\minus\hat{\imath}}] \\
  & \leq \int_{\mathbf{W}^{\minus\hat{\imath}}}{ \Big( 2^b \hspace*{-0.5em}\sum_{l=0}^{\rho(\set{P}_{\hat{\imath}})-1} \binom{K_{\hat{\imath}}}{l} \epsilon_{\hat{\imath}}^l (1-\epsilon_{\hat{\imath}})^{K_{\hat{\imath}}-l} \Big) \mathbb{P}^{K_{\minus\hat{\imath}}}[ d\, \mathbf{W}^{\minus\hat{\imath}}]} \\
  & = 2^b \sum_{l=0}^{\rho(\set{P}_{\hat{\imath}})-1} \binom{K_{\hat{\imath}}}{l} \epsilon_{\hat{\imath}}^l (1-\epsilon_{\hat{\imath}})^{K_{\hat{\imath}}-l} \,. \qed
\end{align*}
\end{pf}

\def\Elproofname{\protect{PROOF of \refthm{ccp:sec2:thm1}.}}
\begin{pf} \label{ccp:pf:sampling_thm}
This proof uses the Boole-Bonferroni inequalities \cite{prekopa1988} as in \cite{kariotoglou2016}. Consider $K_i$ satisfying \eqref{ccp:eqn:implicit} for $i\in\set{I}$. By \reflem{ccp:lemma:sampling_lemma} we have $\mathbb{P}^K [V_i^\star[\mathbb{K}] > \epsilon_i ] \leq \beta_i$ for $i\in\set{I}$, and by subadditivity of $\mathbb{P}^K$ it follows that
\begin{equation*}
   \mathbb{P}^K \big[ \bigcap_{i=1}^{P}\nolimits \{V_i^\star[\mathbb{K}] \leq \epsilon_i \} \big] \geq 1 - \sum_{i=1}^{P} \mathbb{P}^K [ V_i^\star[\mathbb{K}] > \epsilon_i ]  \geq 1 -\mathds{1}^\transp \greekbf{\beta} \geq 1 - \beta\,.
\end{equation*}
It remains to show that
\begin{equation} \label{ccp:eq:sce2:thm2:eq1}
  \mathbb{P}^K [V^\star[\mathbb{K}] \leq \epsilon ] \geq  \mathbb{P}^K \big[ \bigcap_{i=1}^{P}\nolimits \{V_i^\star[\mathbb{K}] \leq \epsilon_i \} \big]\,.
\end{equation}
By \refdef{ccp:def:constraintViolationProbability} we have that
\begin{align*}
  V^\star[\mathbb{K}] &=  \mathbb{P}\{ \mathbf{w} \in \setbf{W} \:\vert\: \exists i \in \mathcal{I}, j \in \set{P}_i \suchthat f_j(\mathbf{x}^\star[\mathbb{K}], \mathbf{y}^\star[\mathbb{K}], \mathbf{w}) > 0  \} \\
  & \leq \sum_{i=1}^{P} \mathbb{P} \{ \mathbf{w} \in \setbf{W} \:\vert\: \exists j \in \set{P}_i \suchthat f_j(\mathbf{x}^\star[\mathbb{K}], \mathbf{y}^\star[\mathbb{K}], \mathbf{w}) > 0 \}\\
  & = \sum_{i=1}^{P} V_i^\star[\mathbb{K}]\,,
\end{align*}
where the inequality follows by subadditivity of $\mathbb{P}$. 
Feasibility, i.e., $V_i^\star[\mathbb{K}]\leq \epsilon_i$ for all $i \in \mathcal{I}$, implies that $V^\star[\mathbb{K}] \leq \mathds{1}^\transp \greekbf{\epsilon} \leq \epsilon$. This implies \eqref{ccp:eq:sce2:thm2:eq1} and yields the result. \hfill\(\qed\)
\end{pf}
   
\end{document}